\input amstex
\documentstyle{amsppt}
\voffset =-3pc



\magnification=\magstep1
\font\chapf=cmbx10 scaled \magstep2

\NoBlackBoxes

\define\rr0{\operatorname{RR\;0}}

\define\conv{\operatorname{conv}}
\define\tsr{\operatorname{tsr}}
\define\RR{\operatorname{RR}}
\define\hull{\operatorname{hull}}

\define\dist{\operatorname{dist}}

\define\topdim{\operatorname{top\;dim}}
\define\locdim{\operatorname{loc\;dim}}
\define\sr1{\operatorname{tsr\;1}}
\define\sa{\operatorname{sa}}
\redefine\sp{\operatorname{sp}}

\redefine\phi{\varphi}
\redefine\epsilon{\varepsilon}
\define\plim{\varprojlim}

\define\wtX{\widetilde{X}}
\define\wtB{\widetilde{B}}
\define\wtC{\widetilde{C}}

\define\wtA{\widetilde{A}}

\define\QA{A^{-1}_q}
\define\bb1{\text{\bf 1}} 
\define\hA{A^\vee}
\document

\topmatter

\title
{\chapf Limits and C*-Algebras of low \\
rank or  dimension}
\endtitle
\rightheadtext{Limits and $C^*-$algebras}
\author
Lawrence G. Brown \& Gert K. Pedersen
\endauthor
\date{Dedicated to the memory of Gert K. Pedersen}\enddate

\abstract{We explore various limit constructions for $C^*$-algebras, such 
as composition series and inverse limits, in relation to the notions of 
real rank, stable rank, and extremal richness.  We also consider 
extensions and pullbacks.  We identify some conditions under which the 
constructions preserve low rank for the $C^*$-algebras or their 
multiplier algebras.  We also discuss the version of topological 
dimension theory appropriate for primitive ideal spaces of 
\,\,$C^*$-algebras and provide an analogue for rank of the countable sum 
theorem of dimension theory.  As an illustration of how the main 
results can be applied, we show that a $CCR$ algebra has stable rank one 
if and only if it has topological dimension zero or one, and we 
characterize those $\sigma$-unital $CCR$ algebras whose multiplier 
algebras have stable rank one or extremal richness. (The real rank zero 
case was already known.)}
\endabstract
\endtopmatter
\document

\vskip2truecm

\subhead{1. INTRODUCTION}\endsubhead

\bigskip

The concept of dimension for a topological space $\Cal X$ 
originates in the basic fact that manifolds are locally 
homeomorphic to euclidean spaces, which have an obvious 
linear dimension. In the more abstract version given by 
\v Cech's covering dimension of a normal space $\Cal X$, the 
dimension gives conditions under which 
certain functions extend and certain cohomology groups 
vanish.

Regarding a $C^*-$algebra $A$ as the non-commutative 
analogue of $C(\Cal X)$ (or $C_0(\Cal X)$) for a compact (or just 
locally compact) Hausdorff space $\Cal X$, it is natural to 
try to extend the notion of topological dimension of 
$\Cal X$ to the analogous setting. The more so as the
covering dimension of $\Cal X$ is easily characterized in 
terms of elements in $C(\Cal X)$. In \cite{{\bf 38}} Rieffel 
defined the {\it (topological) stable rank}, $\tsr(A)$, 
of an arbitrary $C^*-$algebra $A$, using concepts from 
dimension theory:  If $A$ is unital, $\tsr(A)$ is the smallest
$d$ in $\Bbb N$ (or $\infty$) such that unimodular $d-$tuples
(namely, those which are left invertible as $d\times 1$ matrices)
are dense in $A^d$. Shortly after, the stable rank was
identified with the Bass stable rank of $A$, 
\cite{{\bf 15}}, which
is a purely algebraic concept.
In particular, by an earlier result of Vaserstein, 
\cite {{\bf 43}}, we have $\tsr (C_0(\Cal X))= 
[\tsize{\frac 12}\dim (\Cal X \cup \{\infty\})]+1$, the factor $\tsize{\frac 12}$ 
arising from the use of complex scalars in $C_0(\Cal X)$.

The {\it real rank} of a $C^*-$algebra was introduced in 
\cite{{\bf 5}} as an alternative to Rieffel's stable rank. 
Formally the only difference is that self-adjoint elements 
replace the general elements in Rieffel's definition, but 
this has unexpected consequences, especially for small 
values of the rank. In general one has $\RR(A)\le 2\tsr (A)-1$, 
and -- pleasing for the eye -- $\RR(C_0(\Cal X))=\dim (\Cal X \cup \{\infty\})$. 
However, in the lowest possible cases, $\tsr(A)=1$ and 
$\RR(A)=0$, the two notions are independent: one may
be satisfied without the other.

One of the real surprises is the symmetry with which stable rank one
and real rank zero sometimes interact with the two $K$-groups for a unital
$C^*-$algebra $A$:
If $I$ is a closed ideal in $A$ and $\tsr(A)=1$, the
natural map $K_0(I)\to K_0(A)$ is injective, whereas the map
$K_1(I)\to K_1(A)$ is injective if $\RR (A)=0$.  Also, the natural
map from Murray-von Neumann equivalence classes of projections in
$A$ to $K_0(A)$ is injective if $A$ has stable rank one, whereas
its image generates the whole group if $A$ is of real rank zero.

Recall from \cite{\bf 6} that a unital $C^*-$algebra $A$ is 
{\it extremally rich} if the open set $\QA$ of 
{\it quasi-invertible} elements is dense in $A$. Here $\QA$ 
is defined as $A^{-1}\Cal E(A)A^{-1}$, where $\Cal E(A)$ 
denotes the set of extreme points in the closed unit ball, 
$A_1$, of $A$. Equivalently, cf\. \cite{\bf 7}, $A$ is 
extremally rich if $\conv(\Cal E(A))=A_1$, so that -- 
as a Banach space -- $A$ has the $\lambda$-property, 
cf\. \cite{\bf 34}. If $A=C(\Cal X)$, extremal richness is 
equivalent to  $\dim(\Cal X)\le 1$. In general, extremal 
richness is a generalization of Rieffel's notion of
stable rank one suitable for not necessarily finite 
$C^*-$algebras. Thus every purely infinite, simple 
$C^*-$algebra is extremally rich, as is every von  
Neumann algebra.  

The low ranks, i\.e\. stable rank one, real rank zero, and extremal
richness, have different formal properties from the higher ranks.
In particular, the low ranks are invariant under Rieffel-Morita
equivalence, but $\tsr(A\otimes \Bbb K)=2$ whenever $\tsr(A)>1$ and
$\RR(A\otimes \Bbb K)=1$ whenever $\RR(A)>0$.  In this paper we primarily consider low ranks,
though a few results include higher rank cases.

For general (non-commutative) $C^*-$algebras the relationship between rank
and dimension is an analogy rather than a theorem.  Nevertheless, 
topological dimension theory has some applications, and we provide
a brief treatment in Section 2.  This is largely, but not entirely, just a
matter of using the appropriate results from topology, but, as we
explain, it would be wrong simply to apply covering dimension to
primitive ideal spaces.  Our treatment includes all $C^*-$algebras
with almost Hausdorff primitive ideal spaces, in particular all type
$I$ $C^*-$algebras.  Section 2 also contains a result about rank,
Theorem 2.12, which is analogous to and inspired by the countable
sum theorem of dimension theory.

Section 3 treats $C^*-$inverse limits, mainly those where all maps are
surjective.  We provide a framework for representing certain
multiplier algebras as such inverse limits.  In Theorem 3.9
we prove that low rank is preserved by surjective $C^*-$inverse
limits, and in Theorem 3.13 we prove that real rank zero and
stable rank one are preserved in certain multiplier algebras which
are non-surjective inverse limits.

Section 4 has results on low rank of pullbacks where at least one
of the maps is surjective.  Via the Busby construction this leads to 
results about low ranks of extensions and multiplier algebras of
extensions.  And Section 5 contains the applications to $CCR$ algebras
and concluding remarks and questions.

Determination of which extensions of low rank $C^*-$algebras have low
rank has been a matter of continuing interest to many
mathematicians.  In all three cases $A$ has low rank if and only if
$I$ and $A/I$ do and an additional hypothesis is satisfied.  For 
real rank zero and stable rank one the additional hypothesis is
just the natural lifting condition, but for extremal richness it is
the natural lifting condition plus a technical hypothesis.  In all
three cases it is desirable to identify circumstances in which the
hypotheses can be simplified.  In the extremal richness case this means 
more than merely eliminating the technical hypothesis.  Results
and remarks on this subject occur in 4.4, 4.6-4.8, and 5.11-5.14.

The authors previously announced a paper entitled, ``Extremally rich
ideals in $C^*-$algebras.''  The present paper and \cite {\bf 9}
constitute an expanded version of that paper.

\vskip2truecm

\subhead{2. TOPOLOGICAL DIMENSION AND LOW RANK}\endsubhead

\bigskip

\definition{2.1. Composition Series} (i) Unless expressly 
mentioned, the word {\it ideal} will in this paper designate 
a closed (and therefore *-invariant) ideal in a $C^*-$algebra.
We say that an increasing series $\{I_{\alpha}| 0\le\alpha\le\beta\}$
of ideals of $A$, indexed by a segment of the
ordinal numbers, is a {\it composition
series of ideals} for $A$ if $I_0 = 0$, $I_{\beta} = A$, and 
$I_{\alpha} = (\bigcup_{{\gamma < \alpha}} I_{\gamma})^=$
for limit ordinals $\alpha$.  $\wtA$ denotes the unitization of $A$,
and ``${}^=$'' denotes norm closure.
\medskip

\noindent (ii) R\o rdam shows in \cite{{\bf 39},
  4.1--4.3} that in every $C^*-$algebra $A$ there 
is a largest ideal $I_{\sr1}(A)$ of stable rank one. If 
we define $\alpha(y) = \dist(y,\wtA^{-1})$, then the 
ideal is given by  
$$
\aligned
I_{\sr1}(A) & = \{ x \in A  \mid  \forall y
\in \wtA : \quad \alpha (x+y) = \alpha (y) \} \\ 
 & = \{x \in A  \mid  x+\wtA^{-1} \subset(\wtA^{-1})^=\}\,.
\endaligned 
$$

Similar constructions are possible with respect to ideals 
of real rank zero and of extremal richness, cf\. 
\cite{{\bf 9}, Theorems 2.3 \& 2.16}. Thus if we let 
$\alpha_{\sa}(z)= \dist(z,\wtA^{-1}_{\sa})$ and define 
$$
\aligned
R &= \{x \in A_{\sa} \mid  \forall y\in \wtA_{\sa} :
\quad \alpha_{\sa}(x+y) = \alpha_{\sa}(y) \}\\
  &= \{x \in A_{\sa} \mid x + \wtA^{-1}_{\sa} 
\subset (\wtA^{-1}_{\sa})^=\}\,,
\endaligned
$$
then $I_{\rr0}(A) = R + iR $  is an ideal of 
real rank zero in $A$, and the largest such.
\medskip

\noindent (iii) It may happen, of course, 
that $A/I_{\sr1}(A)$ has a non-zero ideal of stable rank 
one (consider e\.g\. the Toeplitz algebra, $\Cal T$), or that 
$A/I_{\rr0}(A)$ has a non-zero ideal of real rank zero 
(consider e\.g\. a non-trivial extension (of real rank
one) of a stabilized Bunce-Deddens algebra by $\Bbb C$, 
arising from a non-liftable projection in its corona 
algebra). In the general case we therefore obtain a strictly 
increasing series $\{I_{\alpha}| 0\le\alpha\le\beta\}$  
of ideals of $A$, which is a composition series for $I_{\beta}$,
such that $I_{\alpha + 1}/I_{\alpha} = I_{\sr1}(A/I_{\alpha})$
and $I_{\sr1}(A/I_{\beta}) = 0$,
and similarly for the real rank zero case. 
\medskip

\noindent (iv) If a $C^*-$algebra $A$ has a  composition 
series of ideals  such that every 
subquotient $I_{\alpha +1}/I_\alpha$ has stable rank one 
(This implies $I_{\beta} = A$ above.) 
we say that $A$ has {\it generalized stable rank one}. 
Similarly we say that $A$ has {\it generalized real rank 
zero} if it has a composition series such that 
$\RR(I_{\alpha +1}/I_\alpha)= 0$ for all $\alpha$.
\medskip

\noindent(v) If $A$ has generalized stable rank one 
and if we choose $I_{\alpha +1}$ such that 
$I_{\alpha +1}/I_\alpha= I_{\sr1}(A/I_\alpha)$, we 
obtain an {\it essential} composition series, 
i\.e\. $I_{\alpha +1}/I_\alpha$ is an essential 
ideal in $A/I_\alpha$ for all $\alpha$. For if $I$ 
were a non-zero ideal of $A/I_\alpha$ orthogonal to 
$I_{\alpha +1}/I_\alpha$, then there is a first 
ordinal $\mu$ such that $J=(I\cap I_\mu)/I_\alpha 
\ne 0$. Since $\mu$ can not be a limit ordinal we 
see that $J$ embeds as an ideal in $I_\mu/I_{\mu -1}$, 
and thus $\tsr(J)=1$. As $J\cap (I_{\alpha +1}
/I_\alpha) =0$ this contradicts our choice of 
$I_{\alpha +1}$ as the largest ideal such that 
$\tsr(I_{\alpha +1}/I_\alpha)= 1$.
Similar reasoning applies in the real rank zero case.
\medskip

\noindent(vi) It follows from the observations made 
above that a $C^*-$algebra $A$ has generalized 
stable rank one or generalized real rank zero 
precisely when $I_{\sr1}(I/J)\ne 0$ or 
$I_{\rr0}(I/J)\ne 0$, respectively, for every non-zero 
quotient $I/J$ of ideals of $A$.

This means that if $A$ has another composition series
$\{I_\alpha \mid 0\le \alpha\le\beta\}$ (determined by 
other interesting subquotient properties), then we can 
find a composition series such that each of its 
subquotients has stable rank one or real rank zero, 
respectively, and is also a subquotient of one of the 
algebras $I_{\alpha+1}/I_\alpha$.
\enddefinition

\bigskip  

\definition{2.2. The Primitive Ideal Space} (i) Recall from 
\cite{{\bf 16}, pp. 233-241} or \cite{{\bf 13}, \S 3} 
that the set $\hA$ of primitive ideals in a
$C^*-$algebra $A$ is a locally quasi-compact, not necessarily 
Hausdorff topological space with the Jacobson topology, 
defined by the closure operation:
$$
\overline{\Cal S} = 
\hull(\ker(\Cal S))\,,\quad \Cal S\subset \hA\,.
$$
Here $\hull (I)= \{P\in \hA \mid I\subset P\}$, 
and $\ker (\Cal S) = \bigcap_{P\in \Cal S} P$ for 
any ideal $I$ of $A$ and any subset $\Cal S$ 
of $\hA$. We obtain the 
formulae $(A/I)^\vee=\hull (I)$ and $I^\vee= 
\hA\setminus \hull (I)$, together with 
$\ker(\hull (I))=I$. Furthermore, for each $x$ in 
$A$ the norm function $\check x$ on $\hA$
given by  $\check x(P) = \Vert x-P\Vert$, 
$P\in \hA$, is lower semicontinuous, so that each 
set $\{P\in \hA\mid \check x (P)> \epsilon\}$ is 
open; and $\check x$ vanishes at infinity, so that the set 
$\{P\in \hA\mid \check x (P)\ge \epsilon\}$  
is compact for $\epsilon > 0$, cf\. 
\cite{{\bf 31}, 4.4.4}.

\medskip

\noindent(ii) For some purposes, including some 
results in Section 5, it is helpful to 
make direct use of topological dimension theory. 
Since by definition $\tsr(A)=\tsr(\wtA)$ and 
$\RR(A)=\RR(\wtA)$ when $A$ is a non-unital 
$C^*-$algebra, we will use $\dim(\Cal X\cup\{\infty\})$ 
as the basic dimension function for any locally 
compact Hausdorff space $\Cal X$. Here $\dim$ is  
\v Cech's covering dimension and $\Cal X\cup\{\infty\}$ 
is the one-point compactification of $\Cal X$. It 
follows from \cite{{\bf 30}, 3.5.6} that 
$$
\dim(\Cal X\cup\{\infty\})=\sup_{\Cal K\subset \Cal X} \{\dim (\Cal K)\}\,,
$$
where $\Cal K$ is compact, and we see from 
\cite{{\bf 30}, 3.5.3} that $\dim(\Cal X\cup\{\infty\})
=\dim (\Cal X)$ when $\Cal X$ is $\sigma-$compact. The
concept of {\it local dimension} is treated in 
Chapter 5 of \cite{\bf 30}, and it follows from 
standard results that $\locdim(\Cal X)=
\dim(\Cal X\cup\{\infty\})$ whenever $\Cal X$ is locally 
compact and Hausdorff.
\medskip

\noindent(iii) Recall that a subset $\Cal S$ of a topological 
space $\Cal X$ is called {\it locally closed} if $\Cal S=\Cal F\cap \Cal G$ 
for some closed and open subsets $\Cal F$ and $\Cal G$ of $\Cal X$, 
respectively. For a $C^*-$algebra $A$ a locally closed 
subset $\Cal S$ of $\hA$ corresponds to a 
subquotient of the form $I/J$, where $I$ and $J$ 
are ideals of $A$ such that $J\subset I$ 
and $\Cal S=(I/J)^\vee = I^\vee\; \cap\; \hull J$. 
Here $I$ and $J$ are not uniquely determined by $\Cal S$, but
$I/J$ is determined up to canonical isomorphism.

\medskip

\noindent(iv) Recall further that a topological 
space $\Cal X$ is called {\it almost Hausdorff} if every 
non-empty closed subset $\Cal F$ contains a non-empty 
relatively open subset $\Cal F\cap \Cal G$ (so $\Cal F\cap \Cal G$ is 
locally closed in $\Cal X$) which is Hausdorff. If $A$ 
is a $C^*-$algebra of type I then $\hA$ is 
almost Hausdorff since every non-zero quotient 
contains a non-zero ideal with continuous trace, 
cf\. \cite{{\bf 31}, 6.2.11}.

\medskip

\noindent(v) We define the {\it topological dimension}, 
$\topdim(A)$, of a $C^*-$algebra $A$ for which $\hA$ is 
almost Hausdorff by 
$$
\topdim (A)=\sup_{\Cal S\subset\hA} \{\locdim (\Cal S)\} 
=\sup_{\Cal K\subset\hA}\{\dim(\Cal K)\}\,,
$$
where $\Cal S$ is any locally closed Hausdorff subset 
and $\Cal K$ is any locally closed compact Hausdorff 
subset of $\hA$.
\enddefinition

\bigskip 

\example{2.3. Remark} In the simplest case where 
$A$ is unital and $\hA$ is Hausdorff we have
$$
\topdim(A)=\dim(\hA)=\RR(Z(A))
$$
by the Dauns-Hofmann Theorem, \cite{{\bf 31}, 4.4.8}, 
where $Z(A)$ denotes the center of $A$. From 
\cite{{\bf 1}} we then deduce, in the case  where  
$A$ is homogeneous of degree $m$ and the corresponding 
Fell bundle is trivial, so that $A=Z(A)\otimes \Bbb M_m$, 
that we have $\RR(A)\le r$ if and only if $\topdim (A)\le (2m-1)r$. 
In particular we learn that it is 
in general false that $\topdim(A)\le \RR(A)$ -- unless 
$\RR(A)=0$, cf\. Proposition 2.9. 
\endexample

\bigskip 

\proclaim{2.4. Proposition} If $I$ is an ideal of a 
$C^*-$algebra $A$ then  $\hA$ is almost Hausdorff 
if and only if $I^\vee$ and $(A/I)^\vee$ 
are both almost Hausdorff, and in that case
$$
\topdim(A)=\max\{\topdim(I), \topdim(A/I)\}\,.
$$
\endproclaim

\demo{Proof} Every open or closed subset of an almost 
Hausdorff space is evidently almost Hausdorff. 
Moreover, any subset of an open or a closed set which 
is relatively locally closed and compact Hausdorff is also 
locally closed and compact Hausdorff in the global 
sense. This proves the first part of the the proposition 
and shows that $\topdim(A)$ majorizes the 
other two. The reverse inequality follows from 
\cite{{\bf 30}, 3.5.6}. \hfill$\square$
\enddemo

\proclaim{2.5. Proposition} If $\{I_\alpha\mid 
0\le \alpha\le\beta\}$ is a composition series for 
a $C^*-$algebra $A$ then $\hA$ is 
almost Hausdorff if and only if $I_{\alpha +1}/I_{\alpha}$
is almost Hausdorff for each $\alpha < \beta$, and if this is so, then
$$
\topdim(A)=\sup_{\alpha < \beta}  \{\topdim(I_{\alpha +1}/I_\alpha)\}\,.
$$
\endproclaim

\demo{Proof} Assume that for some ordinal $\lambda$ 
we have proved that 
$$
\topdim(I_\mu)=\sup_{\alpha < \mu}  \{\topdim(I_{\alpha +1}/I_\alpha)\}
$$
for all $\mu <\lambda$. If $\lambda$ is a limit 
ordinal then $I_\lambda=\left(\bigcup_{\mu<\lambda} 
I_\mu\right)^=$. Since every compact subset of 
$I_\lambda^\vee$ is contained in 
some $I_\mu^\vee$ we conclude that 
$$
\topdim(I_\lambda)=\sup_{\mu<\lambda} \{\topdim(I_\mu)\}
=\sup_{\alpha <\lambda} \{\topdim(I_{\alpha +1}/I_\alpha)\}\,.
$$
If $\lambda$ is not a limit ordinal, i\.e\. $\lambda
=\mu +1$ for some $\mu<\lambda$, then again 
$$
\aligned
\topdim(I_\lambda)&=\max\{ \topdim(I_\lambda/I_\mu), 
\sup_{\alpha<\mu} \{\topdim(I_{\alpha +1}/I_\alpha) \}\}\\
&=\sup_{\alpha<\lambda}  \{\topdim(I_{\alpha +1}/I_\alpha)\}
\endaligned
$$
by Proposition 2.4. The argument can now be completed 
by transfinite induction. \hfill $\square$
\enddemo

\proclaim{2.6. Proposition} If $A$ is a $C^*-$algebra 
such that $\hA$ is almost Hausdorff and if 
$\hA=\bigcup_n \Cal S_n$ where each $\Cal S_n$ is locally
closed, then $\topdim(A)=\sup_n \{\topdim(A_n)\}$, where 
$A_n$ is the subquotient of $A$ with $A_n^\vee=\Cal S_n$.
\endproclaim

\demo{Proof} It follows from Proposition 2.4 that 
$\topdim(A_n)\le \topdim(A)$ for every $n$.

To prove the reverse inequality we shall use 
transfinite induction to construct a composition 
series $\{I_\alpha\mid 0\le\alpha\le\beta\}$ for $A$
such that each locally closed subset 
$(I_{\alpha +1}/I_\alpha)^\vee$ is contained in some $\Cal S_n$. 

Assume that for some ordinal $\lambda$ we have defined 
these ideals for all $\mu<\lambda$. If $\lambda$ is a 
limit ordinal we just put $I_\lambda = 
(\bigcup_{\mu<\lambda} I_\mu)^=$. If $\lambda$ is not 
a limit ordinal, i\.e\. $\lambda=\mu +1$ for some 
$\mu$ such that $I_{\mu} \ne A$, we observe that 
$$
\hull(I_\mu)=(A/I_\mu)^\vee=\bigcup_n  \Cal S_n\cap\hull(I_\mu)\,.
$$
Since $(A/I_\mu)^\vee$ is a Baire space there is 
an $n$ such that $(\Cal S_n\cap\hull(I_\mu))^-$ has non-empty 
interior (relative to $\hull (I_\mu)$). As any locally 
closed subset $\Cal S$ is a dense, relatively open subset
of $\overline{\Cal S}$, this implies that $\Cal S_n\cap \hull(I_\mu)$ 
has a relative interior $\Cal G \ne \emptyset$. 
Thus we may define $I_\lambda$ such that 
$(I_\lambda/I_\mu)^\vee=\Cal G \subset \Cal S_n$.

Since $I_\lambda$ is strictly larger than $I_\mu$, the 
inductive process must eventually terminate with 
$I_\beta=A$ for some ordinal $\beta$, giving us the 
desired composition series. The result now follows 
from Proposition 2.5. \hfill$\square$ 
\enddemo

\bigskip 

\example{2.7. Remarks} (i) We have not been able to  
locate a precise reference for the topological  
analogues of 2.5 and 2.6, but we note that 2.5 is 
a special case of 2.6 when the composition series is 
countable and that the topological analogue of 2.6 
follows from standard results when $\Cal X$ is
second countable.

\medskip

\noindent(ii) It follows from 2.5 that when  
$\hA$ is almost Hausdorff then 
$\dim(\Cal K)\le\topdim(A)$ for any compact Hausdorff 
subset of $\hA$, regardless of whether $\Cal K$ is 
locally closed or not. We do not know whether
such non-locally closed subsets can actually exist.

\medskip

\noindent(iii) Note that $\topdim(A)$ depends only on
$\hA$, but we are not showing this in the notation because
it is not the same as $\dim(\hA)$.
It would be a mistake if we had simply 
defined $\topdim(A)$ to be $\dim(\hA)$, 
ignoring the fact that 
$\hA$ may not be Hausdorff. To see this let 
$\Cal X$ be a locally compact Hausdorff space and define 
$A=(C_0(\Cal X)\otimes \Bbb K)\,\widetilde{}$. Then
$\hA=\Cal X\cup\{\infty\}$, where $\Cal X$ has the given 
topology; but $\Cal X\cup\{\infty\}$ is not the 
one-point-compactification. In fact, the only
open set containing $\infty$ is the whole of 
$\hA$. Thus $\dim(\hA)=0$, whereas we have 
correctly defined $\topdim(A)$ to be $\locdim(\Cal X)$. 
Note also that $\hA$ is compact, so 
$\locdim(\hA)$, as defined in \cite{{\bf 30}}, 
is zero. The space $\hA$, for $\Cal X=[0, 1]$,
appears in \cite{{\bf 30}, Example 3.6.1}.

The phenomenon above occurs annoyingly often, 
and means that we have to work with non-unital 
$C^*-$algebras in order not to destroy the 
Hausdorff properties of their primitive ideal 
spaces.
\medskip

\noindent (iv) Kirchberg and Winter \cite{\bf 17} have defined the decomposition rank, dr$(A)$,
 for
nuclear $C^*-$algebras $A$ and have presented
strong evidence that it is a non-commutative analogue of topological
dimension.  Unlike $\topdim(A)$, dr$(A)$ is not a property
of $\hA$, and therefore it can give much deeper
information about $A$ than can $\topdim(A)$.  Winter \cite{\bf 44} has
shown that dr$(A) =\topdim(A)$ when $A$ is subhomogeneous, but this does not
hold for all type $I$ $A$, since by \cite{{\bf 17}, Example 4.8}, dr$(\Cal T) = \infty$.
\medskip

\noindent (v)  It can be shown that for type $I$ $A$, $\topdim(A) \le \sup \{
\topdim(A_{\alpha})\}$ when $A= (\bigcup A_{\alpha})^=$ for an upward directed
family $\{A_{\alpha}\}$ of $C^*-$subalgebras, cf\. \cite{{\bf 42}, Axiom X3}.
\medskip

\noindent (vi) Although we have defined the topological dimension
of $A$ only when $\hA$ is almost Hausdorff, the concept can be
extended to arbitrary $A$ in the special case of dimension zero.
Thus we define $\topdim (A)=0$ to mean that $\hA$ has a basis
consisting of compact-open subsets.  This concept has already
been used to good effect by Bratteli and Elliott \cite {\bf 2} and 
recently by Pasnicu and R\o rdam \cite {\bf 29}.  In \cite{{\bf 29}, Corollary 4.4}
it is shown that the analogue of 2.4 holds for the extended
concept in the separable case.  However, the proof depends on the
main theorem of \cite {\bf 29}.  We provide in the next proposition a simple
direct proof that the analogue of 2.4 always holds.  Then it is
routine to show that the new definition for $\topdim (A)=0$ agrees
with that given in 2.2 when $\hA$ is almost Hausdorff and that 2.5 
and 2.6 still hold for the new definition.  We are grateful to M.
R\o rdam for providing a copy of \cite {\bf 29} and for helpful discussions.
\endexample

\bigskip

\proclaim{2.8. Proposition} Let $\Cal X$ be a locally quasi-compact
topological space and $\Cal F$ a closed subset.  Then $\Cal X$ has
a basis of compact-open sets if and only if both $\Cal F$ and 
$\Cal X\setminus \Cal F$ have bases of (relatively) compact-open
sets.
\endproclaim

\demo{Proof} Since one direction is trivial, it is enough to
assume $\Cal F$ and $\Cal X\setminus \Cal F$ have the property
and prove that $\Cal X$ does.  Thus we are given a point $p$ in
$\Cal X$ and an open set $\Cal U$ containing $p$, and we need to find
a compact-open set $\Cal C$ such that $p\in \Cal C$ and $\Cal C\subset \Cal U$.
We assume $p$ is in $\Cal F$, since otherwise the existence of $\Cal C$
is obvious.  Then there is a compact relatively open subset $\Cal K$ of
$\Cal F$ such that $p\in \Cal K$ and $\Cal K\subset \Cal U$.  Let
$\Cal V$ be an open set such that $\Cal V\subset \Cal U$ and 
$\Cal V\cap \Cal F=\Cal K$.  By local quasi-compactness there is a
compact set $\Cal L$ such that $\Cal L\subset \Cal V$ and 
$\Cal K\subset \Cal L^{\circ}$, where $\Cal L^{\circ}$ is the
interior of $\Cal L$.  Then $\Cal L\setminus \Cal L^{\circ}$ is a
compact subset of $\Cal V\setminus \Cal F$.  Thus we can find a
compact-open subset $\Cal C_1$ (open relative to $\Cal X\setminus \Cal F$
and hence open in $\Cal X$) such that 
$\Cal L\setminus \Cal L^{\circ}\subset \Cal C_1\subset \Cal V\setminus \Cal F$.
Then let $\Cal C= \Cal L^{\circ}\cup \Cal C_1=\Cal L\cup \Cal C_1$.
\hfill$\square$
\enddemo

\bigskip 

Recall that $A$ is said to have the {\it ideal property} if every
ideal of $A$ is (ideally) generated by its projections.  This
property was defined by Stevens \cite {\bf 41} and extensively studied by
Pasnicu, cf\. \cite{\bf 27}.  A weaker property is that
$A\otimes \Bbb K$ have the ideal property.

\proclaim{2.9. Proposition} If $A$ is a $C^*-$algebra of 
generalized real rank zero, or more generally if $A$ has a
composition series $\{I_\alpha\mid 0\le \alpha \le \beta\}$
such that $(I_{\alpha +1}/I_\alpha)\otimes \Bbb K$ has the ideal property 
for each $\alpha < \beta$, then $\topdim (A)=0$.
\endproclaim

\demo{Proof} Since $A$ and $A\otimes \Bbb K$ have the same primitive
ideal spaces, and since topological dimension is compatible with
composition series (2.5 and 2.7(vi)), it is sufficient to assume $A$
has the ideal property and prove $\hA$ has a basis of compact-open
sets.  For this we use the fact that for every projection $p$ in
$A$, $\{P\in \hA\mid p\notin P\}$ is a compact-open subset of $\hA$.
Thus the ideal property yields directly the fact that
every open subset of $\hA$ is a union of compact-open sets.
\hfill$\square$
\enddemo

\bigskip 

\example{2.10. Remark} Since the study of ranks is our 
primary concern, we have only introduced the topological 
dimension for $C^*-$algbras as a tool. But it raises some 
natural questions. One of the main theorems in topological 
dimension theory is that if a normal space $\Cal X$ is written 
as $\Cal X=\bigcup \Cal F_n$, where each $\Cal F_n$ is closed, then 
$\dim(\Cal X)=\sup_n \dim(\Cal F_n)$. We proceed to establish an 
analogue of this countable sum theorem. Another will 
appear in Theorem 3.9, cf\. 3.10, 3.11. 
\endexample

\bigskip 

\definition{2.11. Definition} Recall from \cite{{\bf 9}} 
that a $C^*-$algebra $A$ is {\it isometrically rich} if 
the union of the left and right invertible elements of 
$\wtA$ is dense in $\wtA$. Equivalently, cf\. 
\cite {{\bf 9}, Proposition 4.2}, 
$A$ is isometrically rich if it is extremally rich and 
$\Cal E(\wtA)$ consists only of isometries and co-isometries. 
\enddefinition

\bigskip 

\proclaim{2.12. Theorem} Let $(I_n)$ be a sequence 
of ideals in a $C^*-$algebra  $A$ such that 
$\bigcup_{n=1}^\infty \hull (I_n)=\hA$. Then: 
\medskip

\itemitem {\rm (i)}  $\RR (A)=\sup_n\{\RR (A/I_n)\}$.
\medskip

\itemitem {\rm (i$'$)}  If $A$ is $\sigma$-unital and $\RR (M(A/I_n))=0$ for all $n$, then $\RR (M(A))=0$.
\medskip

\itemitem {\rm (ii)} $\tsr (A) = \sup_n \{\tsr (A/I_n)\}$.
\medskip

\itemitem {\rm (ii$'$)}  If $A$ is $\sigma$-unital and $\tsr (M(A/I_n))=1$ for all $n$, then $\tsr (M(A))=1$.
\medskip

\itemitem {\rm (iii)} If each quotient $A/I_n$ is isometrically 
rich and $I_{n+1}\subset I_n$ for all $n$, then $A$ 
is isometrically rich.
\medskip

\itemitem {\rm (iv)} If each quotient $A/I_n$ is extremally 
rich and either $\{I_n\}$ is finite or $\hA$ is Hausdorff, then 
$A$ is extremally rich.
\endproclaim

\demo{Proof} Without loss of generality we may 
assume in cases (i), (ii) and (iii) that 
$A$ is unital. 
In case (i) we assume for some $d\ge 0$ that 
$\RR (A/I_n)\le d$ for all $n$ and take a tuple 
$\underline x$ in $(A_{\sa})^{d+1}$. We then wish 
to approximate $\underline x$ by a unimodular 
self-adjoint $(d+1)-$tuple. In case (ii) we 
assume for some $d\ge 1$ that $\tsr (A/I_n)\le d$ 
for all $n$ and take a tuple $\underline x$ in 
$A^d$. We then seek the same kind of approximation 
by a unimodular $d-$tuple. In cases (iii) and (iv) 
we take an arbitrary element $x$ of $A$ or $\wtA$ 
and wish to approximate it by a one-sided invertible 
element of $A$ or by a general quasi-invertible 
element of $\wtA$. The basic construction is the same 
in all these cases, so we will write it out in case (iii), 
the most difficult, and then indicate the minor changes 
to be made for the others.

Recall that for a quasi-invertible element $y$ in 
$A$ we have defined $m_q(y)=\dist(y, A\setminus A_q^{-1})$, 
cf\. \cite{{\bf 6}, 1.4 \& 1.5}, and that $m_q(y)$ 
also measures the distance from $0$ to the rest 
of the spectrum of $|y|$. Now take $\epsilon >0$ 
and let $x_0=x$. Also, let $\pi_n\colon A @>>> A/I_n$ 
denote the quotient morphism. By induction we will 
construct a sequence $(x_n)$ in $A$, such that 
$\pi_n(x_n)\in (A/I_n)_q^{-1}$ and 
$\Vert x_n - x_{n-1}\Vert < \delta_n$, where 
$$
\delta_n = \min\{ 2^{-n}\epsilon\,, \tfrac12 m_q(\pi_{n-1}(x_{n-1}))\,, 
\tfrac 12\delta_{n-1}\} 
$$
for all $n$ (with the convention that $m_q(\pi_0(x_0))=
\delta_0=1$). Assume that we have defined $x_k$ 
for all natural numbers $k<n$. 
Since $A/I_n$ is extremally 
rich we can find $y$ in $(A/I_n)_q^{-1}$ such that 
$\Vert y-\pi_n(x_{n-1})\Vert < \delta_n$, and we may
then choose $x_n$ in $A$ such that $\pi_n(x_n)=y$ and  
$\Vert x_n-x_{n-1}\Vert < \delta_n$, completing the 
induction step. 

The sequence $(x_n)$ is evidently convergent, so we 
can define $x_\infty =\lim x_n$. The inequalities 
in the construction imply that $\Vert x_\infty - x\Vert <
\epsilon$. Moreover, 
$$
\Vert x_\infty - x_n\Vert < \sum_{k>n} \delta_k 
\le 2\delta_{n+1}\le m_q(\pi_n(x_n))\,,
$$
so $\pi_n(x_\infty)\in (A/I_n)_q^{-1}$ for all $n$. 
By assumption $A/I_n$ is isometrically rich, so 
$\pi_n(x_\infty)$ is either left or right invertible. 
However, since $I_{n+1}\subset I_n$, if 
$\pi_n(x_\infty)$ is not right invertible then 
$\pi_m(x_\infty)$ cannot be right invertible for 
any $m\ge n$. We may therefore assume that  
$\pi_n(x_\infty)$ is, say, left invertible for all 
$n$. Equivalently,  $\pi_n(x_\infty^* x_\infty)$ is 
invertible for all $n$. Since $\hA=
\bigcup \hull (I_n)$ this implies that  
$\rho(x_\infty^* x_\infty)$ is invertible for 
every irreducible representation $(\rho, \Cal H)$ 
of $A$. Therefore $x_\infty^* x_\infty$ is
invertible in $A$, so that $x_\infty$ is left 
invertible in $A$, as desired.

In case (iv), when $\hA$ is Hausdorff, the same construction will produce 
an approximant $x_\infty$ in $\bold1+A \subset \wtA$ such that  
$\rho(x_\infty)$ is either left or right invertible 
for every irreducible representation $(\rho,\Cal H)$ 
of $A$. Thus 
$$
m_q(\rho(x_\infty))=
\max\{m(\rho(x_\infty)), m(\rho(x^*_\infty))\} > 0\,,
$$
where, as usual,
$$
m(\rho(x_\infty))=\max\{\epsilon > 0 \mid [0, \epsilon[ \cap \sp
(\rho(|x_\infty|))=\emptyset\}.
$$
Thus $\min(m(\rho(x_\infty)),1)=1-\Vert (\bold1-\rho(|x_\infty|))_+\Vert$.
Since this function is continuous on $\hA$ and approaches $1$
at $\infty$, there is an $\epsilon > 0$ such 
that $m_q(\rho(x_\infty))\ge\epsilon$ for all 
$(\rho, \Cal H)$, whence $m_q(x_\infty)\ge\epsilon$ 
and $x_\infty \in \wtA^{-1}_q$, cf\. \cite{{\bf 6}, 
Proposition 1.2}.  And of course the completion of case (iv) when
$\{I_n\}$ is finite is obvious.

In cases (i) and (ii) a key fact is that the 
set of (self-adjoint) unimodular tuples in $A^d$ 
is open. The distance of a unimodular tuple 
to the set of non-unimodular tuples will then 
replace the function $m_q$ in the previous argument. 
We therefore obtain a tuple $\underline x_\infty 
=(y_1,\dots, y_d)$ such that $\rho(\sum y^*_ky_k)$ 
is invertible for every irreducible representation 
of $A$, which means that $\sum y^*_ky_k$ is 
invertible in $A$ and $\underline x_\infty$ is 
a unimodular tuple.

The proofs of cases (i$'$) and (ii$'$) are almost identical,
so we write it out only for case (ii$'$), the more difficult.
Thus we are given $x$ in $M(A)$ and seek to approximate $x$ by
an invertible element $y$ of $M(A)$.  By the same basic argument
as above, we can approximate $x$ by $x_\infty$ in $M(A)$ such that
for all $n$ $\overline{\pi}_n(x_{\infty})$ is invertible in 
$M(A/I_n)$.  Here $\overline{\pi}_n\colon M(A) @>>> M(A/I_n)$ is the
natural extension of $\pi_n$, and we use the non-commutative Tietze
extension theorem, \cite {{\bf 32}, Theorem 10}, to deduce that $\overline{\pi}_n$
is surjective.  We shall show that $x_{\infty}$ has a polar
decomposition, $x_{\infty}=u|x_{\infty}|$, where $u$ is unitary
in $M(A)$.  Then let $y=u(|x_{\infty}|+\epsilon \bold1)$.

For each $n$ let $(\rho_n,\Cal H_n)$ be a non-degenerate
representation of $A$ with kernel $I_n$.  Let $(\rho,\Cal H)$ be the 
direct sum representation, and let $\overline{\rho}_n$ and
$\overline{\rho}$ be the extensions of $\rho_n$ and $\rho$ to
$M(A)$.  Then $\overline{\rho}$ is faithful and $\overline{\rho}(M(A))$ 
is the idealizer of $\rho(A)$ in $B(\Cal H)$. Since
$\overline{\rho}_n(x_{\infty})$ is invertible, we can write
$\overline{\rho}_n(x_{\infty})= U_n\overline{\rho}_n(|x_{\infty}|)
$, where $U_n$
is unitary in $B(\Cal H_n)$.  Let $U=\bigoplus_nU_n$ and note that
$U\overline{\rho}(|x_{\infty}|)=\overline{\rho}(x_{\infty})
=\overline{\rho}(|x_{\infty}^*|)U$.  It is sufficient to show
that $U$ idealizes $\rho(A)$.  Clearly $U\rho(R)\subset \rho(A)$ and
$\rho(L)U\subset \rho(A)$, where $L=(A|x_{\infty}^*|)^=$ and 
$R=(|x_{\infty}|A)^=$, one-sided ideals of $A$.  We claim that $L=A=R$.
To see this, use the corresponding hereditary $C^*-$subalgebras,
$B=(|x_{\infty}^*|A|x_{\infty}^*|)^=$ and
$C=(|x_{\infty}|A|x_{\infty}|)^=$.  If, for example, $R\ne A$, then
$C\ne A$; and hence $\phi_{|C}=0$ for some pure state $\phi$.  But
since $\bigcup_{n=1}^\infty \hull (I_n)=\hA$, $\phi$ factors through
$A/I_n$ for some $n$.  This contradicts the invertibility of
$\overline{\pi}_n(|x_{\infty}|)$.
\hfill $\square$
\enddemo

\bigskip 

\example{2.13. Remarks} (i) It follows from  
2.12 that we have:
\medskip

\noindent (1)\quad $\RR(A/(I\cap J))=\max\{\RR(A/I), \RR(A/J)\}$,
\medskip
 
\noindent (2)\quad $\tsr(A/(I\cap J))=\max\{\tsr(A/I), \tsr(A/J)\}$, and
\medskip

\noindent (3)\quad $A/(I\cap J)$ is extremally rich$ \Leftrightarrow A/I$ and $A/J$ are 
\medskip

\noindent for any pair $I, J$ of ideals in $A$. Since 
$A/(I\cap J)$ is a surjective pullback of $A/I$ and $A/J$, 
(1) and (2) are not new, cf\. 4.1 and 4.2 below.
\medskip

\noindent (ii) Since the extra conditions in cases (iii) and (iv)
are disappointing, we mention a couple of complements.
\medskip

\noindent (4)\quad If $A/I_1$ is isometrically rich and $\tsr(A/I_n)=1$ 
for $n>1$, then $A$ is 

\quad isometrically rich.
\medskip

\noindent (5)\quad If $A/I_1$, $\dots$, $A/I_m$ are isometrically rich,
$\tsr(A/I_n)=1$ for $n>m$, and if 

\quad $\hull (I_1)$, $\dots$, 
$\hull (I_m)$ are mutually disjoint, then $A$ is extremally rich.
\medskip

\noindent If we replace the sequence $(I_n)$ by 
$(\bigcap_{k=1}^n I_k)$, it is easy to deduce (4) from case (iii).
The deduction of (5) is also elementary, though not so
immediate.
\medskip

\noindent (iii) We show in 4.10 that if $A/I_1$ is only extremally
rich in (4), then $A$ need not be extremally rich, and that the
disjointness hypothesis cannot be omitted from (5).  Thus the extra
conditions in cases (iii) and (iv) of 2.12 cannot be omitted.  We
also show in 4.10 that (iii$'$) and (iv$'$), the analogues of (iii)
and (iv) for multiplier algebras, are false.  It can also be shown
that (i$'$) and (ii$'$) would be false for real ranks $>0$ or
stable ranks $>1$.
\medskip

\noindent (iv) Part (i) gives in principle a new proof of
the topological  countable sum theorem for compact 
Hausdorff spaces. However, the standard proof as found in 
\cite{{\bf 30}, Theorem 2.5} is also of a function-theoretic 
nature. But it uses the fact that $\dim \Cal X \le n$ if and only if 
each unimodular $(n+1)-$tuple in any quotient of $C(\Cal X)_{\sa}$ lifts 
to a unimodular tuple in $C(\Cal X)_{\sa}$. Since general 
$C^*-$algebras may not have very many ideals this definition 
does not generalize, and the commutative proof cannot be used 
as it stands. 
\endexample

\vskip2truecm

\subhead{3. LOW RANK OF INVERSE LIMITS}\endsubhead

\bigskip

\definition{3.1. Inverse Limits} (i) If $\{A_i\}$ is a 
family of $C^*-$algebras indexed by a directed set $\Bbb I$, 
and
if $\{\phi_{ij}\mid i,j\in \Bbb I,\;i>j\}$ is a family of morphisms 
$\phi_{ij}\colon A_i\to A_j$ satisfying the coherence relations 
$\phi_{jk}\circ\phi_{ij}=\phi_{ik}$ for all $i>j>k$ in $\Bbb I$, we
define the {\it $C^*-$inverse limit} as the $C^*-$algebra $\plim A_i$
of bounded strings $x=(x_i)$ in $\prod A_i$ such that
$\phi_{ij}(x_i)=x_j$ for all $i>j$. If $\rho_i\colon A\to A_i,\, i\in
\Bbb I$, is a family of morphisms which is coherent with respect to
$(\phi_{ij})$, i\.e\.  $\phi_{ij}\circ \rho_i = \rho_j$ for all $i > j$,
there is a unique morphism  $\rho\colon A\to \plim A_i$ given by
$\rho(x) = (\rho_i(x))$. This universal property provides an
alternative definition of $\plim A_i$.  

We shall here be exclusively interested in the case where  $\Bbb
I=\Bbb N$. If each morphism $\phi_n = \phi_{n,n-1}$ is surjective, 
we shall refer to $\plim A_n$ as the {\it surjective $C^*-$inverse
  limit} of the $A_n$'s.  For the rest of this section we shall assume
$\plim A_n$ denotes a surjective $C^*-$inverse limit unless we
explicitly say otherwise.

In stark contrast to the direct limit, the inverse limit of $C^*-$algebras is
practically absent from the theory. The reason is that it tends to be
unmanageably large. To circumvent this difficulty Phillips considered
in \cite{{\bf 36}} and \cite{{\bf 37}} the category of
pro$-C^*-$algebras in which full inverse limits (containing unbounded
strings in $\prod A_i$) are allowed, but a much weaker topology 
is used. Roughly speaking, this
is the non-commutative analogue of passing from the category of
compact spaces (where an infinite topological union need not be
in the category) to the category of normal spaces (where this process is 
allowed, cf\. \cite{{\bf 30}, 1.4.3}). In fact, if $\Cal X$ is the
direct limit (in the category of topological spaces) of a directed
family $(\Cal X_i)$ of compact Hausdorff spaces then $\plim C(\Cal X_i) =
C(\beta \Cal X)$, where $\beta \Cal X$ denotes the Stone-\v Cech compactification
of $\Cal X$ (so that $\beta \Cal X$ is the direct limit of $(\Cal X_i)$ in the
category of compact topological spaces). We shall here turn the usual
disadvantage of inverse limits into an advantage, describing a rather
general method of writing the multiplier algebra $M(A)$ of some
$C^*-$algebras $A$ as a surjective $C^*-$inverse limit of quotients of $A$. 

\medskip

\noindent (ii) For every $m$ we define the surjective morphism 
$\pi_m\colon \plim A_n
\to A_m$ by evaluating an element $x=(x_n)$ of $\plim A_n$ at
$m$. Note that $\phi_n\circ \pi_n = \pi_{n-1}$ for $n>1$. For each
$C^*-$subalgebra $B$ of $\plim A_n$ we therefore obtain a 
sequence of $C^*-$subalgebras $B_n =\pi_n(B)$ 
such that $\phi_n(B_n)=B_{n-1}$ for 
$n>1$. Conversely, and more to the point, given such a coherent sequence
$(B_n)$, there is a natural embedding of $\plim B_n$ as a
$C^*-$subalgebra $\overline B$ of $\plim A_n$. Evidently $B \subset
\overline B$, and in general the inclusion is strict. It is
straightforward to verify that if $B$ is an ideal or a
hereditary $C^*-$subalgebra of $\plim A_n$, then this is also the case
for every $B_n$ in $A_n$. 
Conversely, if a  coherent sequence $(B_n)$
consists of ideals or hereditary
$C^*-$subalgebras, then $\overline B$ will have the same property in
$\plim A_n$. With the obvious modifications this holds even for
one-sided ideals.   
\enddefinition

\bigskip

\proclaim{3.2. Theorem} If $\plim A_n$ is the surjective
$C^*-$inverse limit of a sequence $(A_n)$ of $\sigma-$unital 
$C^*-$algebras, and if $A_o$ is an ideal of $\plim A_n$ such 
that $\pi_m (A_o) = A_m$ for every $m$, then 
$M(A_o)=\plim M(A_n)$, a surjective $C^*-$inverse limit.
\endproclaim 

\demo{Proof} Put $A=\plim A_n$ and $M=\plim M(A_n)$.
We then claim that there is a commutative diagram
$$
\CD
A_o @>{\iota_o}>> A @>\pi_n>> A_n @>{\phi_{n}}>>          A_{n-1}\\
@VV{\iota}V  @VV{\rho}V   @VV{\iota_{n}}V      @VV{\iota_{n-1}}V \\
M(A_o) @. M @>{\overline{\pi}_n}>> M(A_{n}) 
@>{\overline\phi_n}>> M(A_{n-1})
\endCD
$$
Here $\iota$ and $\iota_k$ for $k\ge 0$ are the natural embeddings,
and $\overline{\phi}_n$ is the surjective morphism obtained from 
\cite{{\bf 32}, Theorem 10}. It follows that $M$ is a surjective 
inverse limit, and the 
coordinate evaluations $\overline{\pi}_n$ are therefore also 
surjective. 

Since the rightmost square of the diagram is commutative, we can 
define the morphism $\rho$ by $(x_n) \mapsto (\iota_n(x_n))$
for every string $x=(x_n)$ in $A$, and we note that 
$\overline{\pi}_n\circ\rho=\iota_n\circ\pi_n$ by this definition. 
Evidently $\rho$ is injective and
$\rho(A)$ is an essential ideal of $M$.

We claim that $A_o$ is essential in $A$. For if $xA_o=0$ for 
some $x$ in $A$, then $\pi_n(x)A_n=0$ for every $n$ by our 
assumption on $A_o$, whence $\pi_n(x)=0$, and therefore $x=0$. 
It follows 
that $\rho(\iota_o(A_o))$ is an essential ideal in $M$. By the 
universal property of multiplier algebras there is therefore 
an injective morphism $\phi\colon M @>>> M(A_o)$ such that 
$\phi\circ\rho\circ\iota_o=\iota$.

Each surjective morphism $\pi_n\circ\iota_o$ extends 
uniquely to a (not yet claimed surjective) morphism 
$\psi_n\colon M(A_o) @>>> M(A_n)$ . Since $\phi_n\circ\pi_n\circ\iota_o
=\pi_{n-1}\circ\iota_o$, 
then also $\overline{\phi_n}\circ\psi_n=\psi_{n-1}$ for 
all $n$. By the universal property of inverse limits this means 
that we have a unique morphism 
$\psi\colon M(A_o) @>>> M$ such that 
$\overline{\pi}_n\circ\psi=\psi_n$ for all $n$. It follows that 
$$
\overline{\pi}_n\circ\rho\circ\iota_o=
\iota_n\circ\pi_n\circ\iota_o=
\psi_n\circ\iota=\overline{\pi}_n\circ\psi\circ\iota
$$
for all $n$, which implies that 
$\psi\circ\iota=\rho\circ\iota_o$.

Combining these results we find that 
$$
(\phi\circ\psi)\circ\iota=\phi\circ\rho\circ\iota_o
=\iota \quad\text{and}
\quad (\psi\circ\phi)\circ\rho\circ\iota_o
=\psi\circ\iota=\rho\circ\iota_o\,. 
$$
Since $\iota(A_o)$ is an essential ideal in $M(A_o)$ and 
$\rho(\iota_o(A_o))$ is an essential ideal in $M$ these equations 
imply that $\phi$ and $\psi$ are the inverses of one another,
and we have our natural isomorphism. \hfill$\square$
\enddemo

\bigskip 

\proclaim{3.3. Corollary} If $(A_n)$ is a sequence of 
$\sigma-$unital $C^*-$algebras with surjective 
morphisms $\phi_n\colon A_n @>>> A_{n-1}$, and if  
$\overline\phi_n\colon M(A_n)@>>> M(A_{n-1})$ denote the 
unique (surjective) extensions of the $\phi_n$'s then 
$$
M(\plim A_n) = \plim M(A_n)\,.
$$
\hfill $\square$
\endproclaim

\definition {3.4. Constant Ideals} We say that an ideal $I$ 
in $\plim A_n$ is $m-${\it constant} if
$I\cap \ker \pi_m = 0$. Equivalently, $I\subset 
(\ker\pi_m)^\perp$. 
Since $(\ker \pi_n)$ is a decreasing sequence,
$I\cap\ker\pi_n=0$ for all $n\ge m$.
Since $\phi_n\circ\pi_n=\pi_{n-1}$ we see that 
$\ker\phi_n \subset \pi_n(\ker\pi_m)$ for $n>m$. If 
therefore $I_n = \pi_n(I)$ denotes the associated sequence 
of ideals in $A_n$, then $I_n\cap\ker\phi_n=0$ for $n>m$.
Thus $I$ is isomorphic to $I_m$ and $I_n$ is isomorphic to 
$I_m$ for all $n\ge m$. In particular, $\plim I_n=I$.
Conversely, if $(I_n)$ is a sequence of ideals in $(A_n)$ 
such that $\phi_n(I_{n})=I_{n-1}$ for $n>1$ and 
$I_n\cap\ker\phi_n=0$ for all $n>m$ for some $m$ then 
$I=\plim I_n$ will be an $m-$constant ideal in 
$\plim A_n$.

If $I$ is an $n-$constant and $J$ an $m-$constant ideal
with $n\le m$, then $I+J\subset(\ker\pi_m)^\perp$ since 
$\ker\pi_m\subset\ker\pi_n$, so $I+J$ is an $m-$constant 
ideal. Since $(\ker\pi_m)^\perp$ is the largest $m-$constant 
ideal, it follows that $A_c=(\bigcup (\ker\pi_m)^\perp)^=$ is
equal to the sum of all constant ideals, and we shall refer 
to it as the {\it quasi-constant ideal} of $\plim A_n$.
\enddefinition

The motivating example for considering constant and 
quasi-constant ideals arises from the Stone-\v Cech 
compactification. If $\Cal X$ is a locally compact Hausdorff 
space then $C_b(\Cal X)$ is always a $C^*-$inverse limit. In 
the case where $\Cal X$ is also $\sigma-$compact we can write 
$\Cal X=\bigcup \Cal X_n$, where each $\Cal X_n$ is compact 
and $\Cal X_n \subset (\Cal X_{n+1})^\circ$. Put $A_n =C(\Cal X_n)$ 
and let $\phi_n(f)= f| \Cal X_n $ for each $f$ in $C(\Cal X_{n+1})$. 
Then $C_b(\Cal X) = \plim A_n$. The large constant ideals will 
be of the form $(\ker\pi_m)^\perp = C_0((\Cal X_m)^\circ)$, so 
the quasi-constant ideal of $C_b(\Cal X)$ can be identified with 
$C_0(\Cal X)$.

Theorem 3.2 provides an immediate generalization of this 
construction:

\bigskip 

\proclaim{3.5. Corollary} If $A = \plim A_n$ is the 
surjective $C^*-$inverse limit of a sequence of $\sigma-$unital 
$C^*-$algebras $(A_n)$, such that the quasi-constant ideal 
$ A_c$ of $A$ satisfies $\pi_m (A_c) = A_m$ for every $m$, 
then $M(A_c) =\plim M(A_n)$. \hfill $\square$
\endproclaim

\bigskip 

\proclaim{3.6. Theorem} Let $(J_n)$ and $(I_n)$ be two 
sequences of ideals in a $C^*-$algebra $A$, one increasing, 
the other decreasing, but such that $I_n\cap J_n=0$ for all 
$n$. If $A_o = (\bigcup J_n)^=$ is essential in $A$ and 
$A_o + I_n = A$ for every $n$, then with $A_n = A/I_n$ 
and $\phi_n\colon A_{n}\to A_{n-1}$ the natural morphisms 
we have an embedding $\rho\colon A @>>> \plim A_n$ such 
that $\rho(A_o)$ is an ideal. If also each $A/I_n$ is
$\sigma$-unital, then $M(A_o)=\plim M(A_n)$.
\endproclaim

\demo{Proof} If $x\in \bigcap I_n$ then it annihilates $J_n$ 
for every $n$, whence $x\in A_o ^\perp$. But then $x=0$ since 
$A_o$ is essential. Thus our assumptions imply that 
$\bigcap I_n = 0$. The quotient morphisms $\rho_n\colon 
A\to A_n$ satisfy $\phi_n\circ\rho_n=\rho_{n-1}$ for all $n$, 
and therefore define a morphism $\rho\colon A\to \plim A_n$. 
Since $\ker\rho =\bigcap I_n = 0$, this is an embedding. 
Observe that $\rho_n(J_m)=(J_m+I_n)/I_n$ is an ideal in 
$A_n$ for every $n$ and $m$. Since $\ker\phi_n=I_{n-1}/I_n$ 
we see moreover that $\rho_n(J_m)\cap\ker\phi_n=0$ for $n>m$. 
Consequently $(\rho_n(J_m))$ is a coherent sequence of ideals 
in $(A_n)$, all isomorphic for $n\ge m$, thus giving rise to 
the $m-$constant ideal $\rho(J_m)$ in $\plim A_n$. It follows 
from this that $\rho(A_o)$ is an ideal in $\plim A_n$, 
isomorphic to $A_o$ (and contained in the quasi-constant ideal 
of $\plim A_n$). Since by assumption
$$
\pi_n(\rho(A_o))=\rho_n(A_o)=(A_o+I_n)/I_n=A_n\,,
$$
it follows from Theorem 3.2 that $M(A_o)=\plim M(A_n)$.
\hfill $\square$
\enddemo

\bigskip

\example{3.7. Remarks} (i) If $A=\plim A_n$ where each $A_n$ is
separable, it can be shown that $A_c$ is the largest separable ideal
of $A$, a result analogous to the result from \cite {\bf 3} that every 
separable $C^*-$algebra $B$ is the largest separable ideal of
$M(B)$.  With the help of these facts it can be shown that if
$M(A_0)=\plim M(A_n)$ for a separable ideal $A_0$ of $A$,
then $A_0=A_c$ and $\pi_m(A_0)=A_m$, as in 3.5.  Also if $(I_n)$ is
a decreasing sequence of ideals of a separable $C^*-$algebra $A$
such that $\bigcap I_n=0$, and if $A_0$ is an ideal of $A$ such
that $M(A_0)$ is identified as above with $\plim M(A/I_n)$, then
$A_0=(\bigcup I_n^\perp )^=$ and $A_0 +I_n =A$, as in 3.6.
These facts provide some justification for our approach in 3.5 and 3.6.

\noindent (ii) In 3.6 we could enlarge $I_n$ to $J_n^\perp$ or
$J_n$ to $I_n^\perp$ and still have the hypotheses.  In the first
case we see that $\plim M(A_n)$ doesn't change.  In the second,
$(\bigcup J_n)^=$ doesn't change, at least in the separable case.
\endexample

\bigskip

\proclaim{3.8. Corollary} Let $(I_n)$ and $(J_n)$ be two 
sequences of ideals in a $C^*-$alge\-bra $A$, one 
decreasing, the other increasing, such that $I_n\perp J_n$ 
for every $n$. If \;$\bigcup J_n$ is dense in $A$ and
each quotient $A_n = A/I_n$ is unital, then with 
$\phi_n\colon A_{n+1}\to A_n$ the natural morphisms we 
have an embedding of $A$ into $\plim A_n$, such that 
$\plim A_n = M(A)$. \hfill $\square$
\endproclaim

\bigskip

\proclaim{3.9. Theorem} Let $A=\plim A_n$ be the surjective 
$C^*-$inverse limit of a sequence of $C^*-$algebras. Then:
\medskip

\itemitem {\rm (i)} If  $\RR(A_n)=0$ for all $n$, then $\RR(A)=0$.
\medskip

\itemitem {\rm (ii)} $\tsr (A)=\sup_n\{\tsr(A_n)\}$.
\medskip

\itemitem {\rm (iii)}  If each $A_n$ is isometrically rich, 
then $A$ is isometrically rich.
\medskip

\itemitem {\rm (iv)} If each $A_n$ is extremally rich, then
$A$ is extremally rich.
\endproclaim

\demo{Proof} We shall use the same basic construction 
as in the proof of Theorem 2.12, relative to the 
surjective morphisms $\pi_n\colon A @>>>A_n$. Assuming, 
as we may, that $A$ is unital we find in cases (i), 
(iii) or (iv) for each $x$ in $A_{\sa}$ or in $A$ and each
$\epsilon >0$ an 
$x_\infty$ in $A_{\sa}$ or in $A$, such that 
$\Vert x_{\infty}-x\Vert <\epsilon$ and 
$\pi_n(x_\infty)\in (A_n)^{-1}_{\sa}$ (in case (i)) and 
$\pi_n(x_\infty) \in (A_n)^{-1}_q$ (in cases (iii) and 
(iv)). Moreover, in case (iii) where each $A_n$ is 
isometrically rich, we see from the connecting morphisms 
$\phi_n\colon A_n @>>> A_{n-1}$ that either all 
$\pi_n(x_\infty)$ are left invertible or all are right 
invertible. 

Let $\pi_n(x_\infty)=w_n|\pi_n(x_\infty)|$ be the polar
decomposition, so that $w_n$ is in $\Cal E(A_n)$. Since 
$w_n$ is unique and $|\pi_n(x_\infty)|=\pi_n(|x_\infty|)$
it follows that $\phi_n(w_n)=w_{n-1}$ for $n>1$, so that 
$w=(w_n)$ is in $\Cal E(A)$. Moreover, if all $\pi_n(x_\infty)$ 
are self-adjoint invertibles, then every $w_n$ is a symmetry, 
so $w$ is a symmetry, and if all $\pi_n(x_\infty)$ are, 
say, left invertible, then each $w_n$ is an isometry, 
so $w$ is an isometry. Now put $y=w(|x_\infty|+\epsilon \bold1)$. 
Then $y\in A^{-1}_q$ (and $y$ is self-adjoint invertible 
if $w$ is a symmetry, whereas $y$ is left invertible if $w$ is 
an isometry) with $\Vert y-x\Vert \le 2\epsilon$.

In the remaining case (ii) we are given  a tuple 
$\underline x$, and the approximant $\underline x_\infty$ 
is also a tuple. Thus $\underline x_\infty= (y_1,\dots, y_d)$ 
is in $A^d$,  and $h_n =\pi_n (\sum y_k^*y_k)$ is invertible 
in $A_n$ for all $n$. Put $\underline w_n =\pi_n(\underline x_\infty)h_n^{-1/2}$. 
Then as above $\phi_n(\underline w_n) =\underline w_{n-1}$
and we approximate $\underline x$ with 
$\underline w((\sum y_k^*y_k)^{1/2} +\epsilon \bold1)$.
\hfill $\square$
\enddemo

\bigskip

\proclaim{3.10. Corollary} Let $(J_n)$ be an 
increasing sequence of ideals in a $C^*-$alge\-bra 
$A$ such that $\bigcup J_n$ is dense in $A$. Assume 
furthermore that each annihilator quotient 
$A/(J_n)^\perp$ is unital. Then:
\medskip

\itemitem {\rm (i)} If $\RR(A) =0$, then $\RR(M(A)) =0$.
\medskip

\itemitem {\rm (ii)} $\tsr(M(A)) =\tsr(A)$.
\medskip

\itemitem {\rm (iii)} If $A$ is isometrically rich, then $M(A)$ is isometrically rich.
\medskip

\itemitem {\rm (iv)} $A$ is extremally rich if and only if
$M(A)$ is extremally rich if and only if each $A/(J_n)^\perp$ is extremally
rich. 
\endproclaim

\demo{Proof} Combine Corollary 3.8 and Theorem 3.9. 
\hfill$\square$
\enddemo

\bigskip

\example{3.11. Remarks$-$Example} 
\noindent (i) The idea in Corollary 3.10 of 
combining properties of ideals and their annihilators 
is found in \cite{{\bf 40}, Proposition 3.15}, which is labeled a technical 
proposition. Viewed as a generalization of writing 
$C_b(\Cal X)$ as a $C^*-$inverse limit $\plim C(\Cal X_n)$, 
cf\. 3.4, the condition seems more intuitive. 
Sheu's result calculates $\tsr(A)$ using weaker hypotheses on
the ideals than those in 3.10.  It helped inspire some of our
results and in turn could be deduced from 2.12(ii).
\medskip

\noindent (ii) It is instructive to realize that these formulae are 
non-commutative analogues of the well-known 
identities $\dim (\beta \Cal X) = \dim (\Cal X)$, valid for
any normal space $\Cal X$, \cite{{\bf 30}, 6.4.3}. By 
contrast, the identities $\tsr (M(A))= \tsr (A)$ 
and $\RR (M(A)) = \RR (A)$ are often false for
non-commutative $C^*-$algebras. A partial ``explanation'' 
might be that $M(A)$ is not always obtainable as a 
$C^*-$inverse limit in the non-commutative case.
\medskip

\noindent (iii) The hypotheses of 3.10 imply
that $\hA =\bigcup \hull (J_n^\perp )^\circ$, a considerably stronger
condition than the one used in 2.12, but cases (iii) and (iv) and case (ii)
for higher ranks don't follow from 2.12.
\medskip

\noindent (iv) The reader may have wondered at 
the asymmetry in the treatment of stable rank 
and real rank in Theorem 3.9 parts (i) and (ii). 
The truth is that we have -- at the moment -- 
no argument to prove that if $\RR(A)\le n$ for 
some unital $C^*-$algebra $A$ and $n>0$, then
for each $\epsilon >0$ there is a $\delta >0$ 
such that for every tuple $(x_1,\dots, x_n)$ in 
$A_{\sa}$ there is a tuple $(y_1,\dots, y_n)$ 
in $A_{\sa}$ with $\sum y_k^2\ge \delta$ and 
$\Vert x_k-y_k\Vert\le\epsilon$ for all $k$. 
In the similar situation for stable ranks we 
can take $\delta$ to be any number less that 
$\epsilon^2$, as we saw in the proof of Theorem 3.9.
This missing information means that the higher 
real ranks of inverse limits and even direct products
cannot be estimated, a fact that seems not to 
be widely known.
\medskip

\noindent (v) We show that the surjectivity hypothesis in Theorem 3.9
cannot be omitted, even in the separable, commutative, unital case.
Note that if $(A_n)$ is a decreasing sequence of $C^*-$subalgebras
of $B$, then $\bigcap A_n$ is the $C^*-$inverse limit of the $A_n$'s.
Let $\Cal X$ be an arbitrary compact metric space and
$f\colon \Cal C\to \Cal X$ a surjective continuous map, where $\Cal C$
is the Cantor set.  Let $\Cal G=\{ (s,t)\in \Cal C\times \Cal C\mid f(s)=f(t)\}$ and
$\Cal D =\{(s_n,t_n)\}$, a countable dense subset of $\Cal G$.  Then if
$B=C(\Cal C)$ and $A_n=\{g\in B\mid g(s_k)=g(t_k), k=1, \dots ,n\}$,
then $\topdim(A_n) =0$ and $\bigcap A_n \cong C(\Cal X)$.
\endexample

Despite this example, we have some positive results about one class of
non-surjective inverse limits:  Let $(I_n)$ be an increasing sequence
of ideals such that $A=(\bigcup I_n)^=$.  Then $M(A)$ is the
$C^*-$inverse limit of the $M(I_n)$'s relative to the restriction maps
$\rho_n\colon M(A)\to M(I_n)$ and $\rho_{n,n-1} \colon M(I_n)\to M(I_{n-1})$.
(But $A$ is the direct limit of the $I_n$'s.)

\bigskip

\proclaim{3.12. Lemma}Let $A=(\bigcup I_n)^=$, where $(I_n)$ is an increasing sequence of
ideals, and let $g\in A_+$.
Then there is an increasing sequence $(g_n)$ such that $g_n\in I_{n+}$ and
$\|g-g_n\|\to 0$.
\endproclaim

\demo{Proof}Let $P_n=\{x\in I_{n+}\colon \|x\| < 1\}$ and $P=\bigcup P_n$.
By \cite {{\bf 31}, 1.4.3} $P_n$ is directed upward and forms an approximate identity of
$I_n$.
Hence $P$ is directed upward and forms an approximate identity of $A$.
Thus we can recursively construct an increasing sequence $(r_j)$ in $P$ such
that $\|g^{1\over 2}(\bold1-r_j) g^{1\over 2}\| < {1\over j}$ and a strictly
increasing sequence $(n_j)$ with $r_j\in P_{n_j}$.
Finally, let $g_n=g^{1\over 2} r_j g^{1\over 2}$ if $n_j\leq n < n_{j+1}$
($g_n=0$ if $n<n_1$).
\hfill $\square$
\enddemo

\proclaim{3.13. Theorem}Assume $A$ is a $\sigma$--unital $C^*$--algebra and $(I_n)$ 
an increasing sequence of ideals such that $A=(\bigcup\limits_n I_n)^=$.
Let $\rho_n\colon M(A)\to M(I_n)$ be the restriction maps.
\medskip

\itemitem {\rm (i)} If $\tsr(A)=1$ and if $p$ and $q$ are projections in $M(A)$ such that
$\rho_n(p)\sim\rho_n(q)$ for all $n$, then $p\sim q$, where $\sim$ denotes
Murray-von Neumann equivalence.
\medskip

\itemitem {\rm (ii)} If $\tsr (M(I_n))=1$ for all $n$, then $\tsr (M(A))=1$.
\medskip

\itemitem {\rm (iii)} If $\RR(M(I_n))=0$ for all $n$, then $\RR(M(A))=0$.
\endproclaim

\demo{Proof}(ii) In the proof we make frequent use of a $C^*$--algebraic operation
which has already occurred in connection with $C^*$--algebras of low rank, but
which has no standard notation.
If $x$ is an element of a $C^*$--algebra $B$ which is faithfully represented on a
Hilbert space, then $x$ has a canonical polar decomposition, $x=v|x|$, where $|x|$
is in $B$ but $v$ need not be.
If $f\colon [0,\infty[\to \Bbb C$ is a continuous function such that $f(0)=0$,
let $x_{[f]}$ denote $vf(|x|)=f(|x^*|)v$.
Then, as is well known, $x_{[f]}$ is in $B$ and is independent of the
representation of $B$.
Also this operation is compatible with morphisms.
For $\delta\geq 0$ let $E_\delta$ and $F_\delta$ be the spectral projections of $|x|$ and $|x^*|$ for the
interval $]\delta,\infty[$, which again need not be in $B$.
Then, when $B$ is unital, R\o rdam \cite{\bf 39} showed that $\tsr(B)=1$ if and only for each $x$ in $B$ and each $\delta>0$ there is a
unitary $u$ in $B$ with $uE_\delta=vE_\delta$ (equivalently $F_\delta u=F_\delta v)$.
In other words $u$ ``extends'' the partial isometry $v E_\delta$.
This property of $u$ is independent of the representation of $B$ and can equivalently be stated:
\medskip

\noindent \,\,\,(6)\quad $x_{[f]}=uf(|x|)$ whenever $f_{|[0,\delta]}=0$, or, still equivalently,
\smallskip

\noindent \,\,\,(7)\quad $xg(|x|)=u |x| g(|x|)$ whenever $g_{|[0,\delta]}=0$.
\medskip

Let $x$ be in $M(A)$ and $0<\delta<1$. Fix a strictly positive element $g$ of $A$, and let $(g_n)$ be as in 3.12.
Define functions $h_0,h,k$, and $\theta_n$ on $[0,\infty[$ by:
\medskip

\noindent \,\,\,(8)\quad $h_{0|[0,{\delta \over 4}]}=0$, $h_{0|[\delta,\infty[}=1$, $h_0$ is linear on $[{\delta \over 4},\delta]$,
\smallskip

\noindent \,\,\,(9)\quad $h_{|[0,{1\over 2}]}=0$, $h_{|[1,\infty[}=1$, $h$ is linear on $[{1\over 2},1]$,
\smallskip

\noindent (10)\quad $k_{|[0,{1\over 4}]}=0$, $k_{|[{1\over 2},\infty[}=1$, $k$ is linear on $[{1\over 4},{1\over 2}]$,
\smallskip

\noindent (11)\quad $\theta_1=k\circ h_0$, and $\theta_{n+1}=k\circ h\circ\theta_n$.
\medskip

Also we shall use the same symbol $\rho_n$ to denote the restriction maps from $M(I_m)$ to $M(I_n)$ for $m\geq n$.

We shall recursively define $w_n,e_n$, and $f_n$ in $M(I_n)$ and $y_n$ in $M(A)$ such that:
\medskip

\noindent (12)\quad $w_n$ is unitary, $0\leq e_n, f_n\leq \bold1$, 
\smallskip

\noindent (13)\quad $w_n$ extends $v_n E_\delta^{(n)}$, referring to a polar decomposition of $\rho_n(x)$,
\smallskip

\noindent (14)\quad $(\rho_{n-1} (w_n)-w_{n-1}) e_{n-1}=0=f_{n-1} (\rho_{n-1} (w_n)-w_{n-1})$ for $n>1$,
\smallskip

\noindent (15)\quad $w_n e_n=f_n w_n$,
\smallskip

\noindent (16)\quad $\|(\bold1-e_n) g_n\|,\ \|g_n(\bold1-f_n)\| < {1\over n}$,
\smallskip

\noindent (17)\quad $\rho_{n-1} (e_n) e_{n-1}=e_{n-1}$ and $\rho_{n-1} (f_n) f_{n-1}=f_{n-1}$ for $n>1$,
\smallskip

\noindent (18)\quad $y_n\in x_{[\theta_n]}+I_n$,
\smallskip

\noindent (19)\quad $\rho_n(y_n)=w_n |\rho_n (y_n)|$,
\smallskip

\noindent (20)\quad $(\bold1-|y_n|) E_\delta=0=F_\delta (\bold1-|y_n^*|)$, and
\smallskip

\noindent (21)\quad $(\bold1-\rho_n (|y_n|)) e_n=0=f_n (\bold1-\rho_n (|y_n^*|))$.
\medskip

To start the construction, choose a unitary $w_1$ which extends $v_1 E_{\delta\over 4}^{(1)}$, establishing (13).
Then let $(r_i)$ be an approximate identity for $I_1$ and let 
$$
e'_i=h(\rho_1 (h_0 (|x|)+(\bold1-h_0 (|x|))^{1\over 2} r_i
(\bold1-h_0 (|x|))^{1\over 2}))\,.
$$
Since $(e'_i)$ converges strictly to $\bold1$ in $M(I_1)$, we may define $e_1=e'_{i_0}$, where $i_0$ is chosen large enough that
$\|(\bold1-e_1) g_1\|$, $\|g_1(\bold1-w_1 e_1 w_1^*)\| < 1$.
Using (15) as the definition of $f_1$, we have (16).
Now let 
$$
c_1=k(\rho_1(h_0 (|x|)+(\bold1-h_0 (|x|))^{1\over 2} r_{i_0} (\bold1-h_0 (|x|))^{1\over 2}))\,,
$$
and $m_1=w_1c_1$.
From the choice of $w_1$ and the fact that $c_1\in\rho_1(\theta_1 (|x|))+I_1$, it follows that $m_1=\rho_1 (x_{[\theta_1]})+z_1$ for some
$z_1$ in $I_1$.
Then let $y_1=x_{[\theta_1]}+z_1$, so that (18) and (19) are clear.
To prove (20) it is enough to show that the relations hold both modulo $I_1$ and
after application of $\rho_1$.
Modulo $I_1$ the relations are obvious from (18) and the fact that $\theta_{1_{|[\delta,\infty]}}=1$.
Since $\rho_1 (y_1)=m_1$, it is clear that $\rho_1 (|y_1|)$ unitizes anything unitized by
$\rho_1(h_0 (|x|))$.
This shows the first part of (20), and the second part follows by conjugation with $w_1$ (using (13)).
The first part of (21) follows from the fact that $kh=h$, and the second again follows by conjugation with $w_1$.

Now let $n>1$ and assume the first $n-1$ steps.
Choose a unitary $w_n$ in $M(I_n)$ which extends $v' E'_{1\over 4}$, referring to a polar decomposition of
$\rho_n(y_{n-1})$.
To verify (13) it suffices to verify it both modulo $I_{n-1}$ and after application of $\rho_{n-1}$.
Then (18) for $n-1$ easily implies the relation modulo $I_{n-1}$.
From (19) for $n-1$ we see:
\medskip

\noindent (22)\quad $\rho_{n-1} (w_n) t=w_{n-1} t$ whenever $\rho_{n-1}( |y_{n-1}|)t=t$, and $t\rho_{n-1} (w_n)=tw_{n-1}$ 

\quad whenever $t\rho_{n-1}(|y_{n-1}^*|)=t$.
\medskip

\noindent In view of (20) and (13) for $n-1$, the proof of (13) is complete.
Next let $(r_i)$ be an approximate identity for $I_n$ and let 
$$
e'_i=h(\rho_n (h(|y_{n-1}|)+(\bold1-h(|y_{n-1}|))^{1\over 2} r_i
(\bold1-h(|y_{n-1}|))^{1\over 2}))\,.
$$
Since $(e'_i)$ converges strictly to $\bold1$ in $M(I_n)$, we may define $e_n=e'_{i_0}$ where $i_0$ is chosen large enough that
$\|(\bold1-e_n) g_n\|$, $\|g_n(\bold1-w_n e_n w_n^*)\| < {1\over n}$.
As in step one, we define $f_n$ by (15) and have (16).
Now we can deduce (14) from (21) for $n-1$, using (22).
Since $\rho_{n-1}(e_n)$ unitizes anything unitized by $\rho_{n-1}(|y_{n-1}|)$, we have the first part of (17).
The second part follows by conjugation, using (14).
Next let 
$$
c_n=k(\rho_n(h(|y_{n-1}|)+(\bold1-h(|y_{n-1}|))^{1\over 2} r_{i_0} (\bold1-h(|y_{n-1}|))^{1\over 2}))\,,
$$  
and $m_n=w_nc_n$.
Since $c_n\in
\rho_n(k(h(|y_{n-1}|)))+I_n$, $m_n\in w_n k(h(\rho_n (|y_{n-1}|)))+I_n$.
Also, the choice of $w_n$ and the fact that $k\circ h$ is supported on $[{1\over 4},\infty)$ imply that $w_n k(h(\rho_n
(|y_{n-1}|)))=\rho_n(y_{n-1})_{[k\circ h]}$.
It follows, using (18) for $n-1$, that $m_n=\rho_n(x_{[\theta_n]})+z_n$ for some $z_n$ in $I_n$.
Then let $y_n=x_{[\theta_n]}+z_n$, so that (18) and (19) are clear.
Finally (20) and (21) are proved as in step one.

Now given the recursion, we construct $w$ essentially as the strict limit of $(w_n)$.
Of course $w_n$ is only in $M(I_n)$, but for each $a$ in $\bigcup I_n$, $w_na$ and $aw_n$ are defined for $n$ sufficiently large; and we
claim that these sequences converge.
In fact from (14) and (17), $f_n\rho_n (w_m)$ and $\rho_n (w_m) e_n$ are constant for $m\geq n$.
In particular $(w_m a)$ is convergent for $a$ in $e_n g_n I_n$.
Since $e_n g_n\to g$ by (16) and the choice of $(g_n)$, we conclude that $(w_m a)$ is convergent for $a$ in $g I_n$; and since $g$
is strictly positive, $(g I_n)^-=I_n$.
The convergence of $(a w_m)$ is proved similarly.
Finally, it is clear that $w$ is unitary and the fact that $w$ extends $v E_\delta$ follows from (13).
\medskip

(iii) It was shown in both \cite {\bf 4} and \cite {\bf 10} that a $C^*$--algebra $B$ has real rank zero if and only if it satisfies an interpolation by projections property.
When $B$ is unital and is faithfully represented on a Hilbert space this property can be stated as follows:\ Let $x$ be in $B_{\sa}$, and for $\delta>0$ let $E_\delta^+$ and $E_\delta^-$ be the spectral projections of $x$ for the intervals $]\delta,\infty[$ and $]-\infty,-\delta[$.
Then there is a projection $p$ in $B$ such that $E_\delta^+\leq p\leq \bold1-E_\delta^-$.
This property can be equivalently stated:
\medskip

\noindent ($7'$)\quad $pg(x)=g(x)$ whenever $g_{|]-\infty,\delta]}=0$, and $pg(x)=0$ whenever $g_{|[-\delta,\infty[}=0$.
\medskip

\noindent Thus the interpolation property is independent of the representation of $B$.
Note that for the canonical polar decomposition, $x=v|x|$, we have that $v^*=v$, $E_\delta=F_\delta=E_\delta^++E^-_\delta$, $vE_\delta=E_\delta^+ - E_\delta^-$, and $v|x|=|x|v$ in the notation of case (ii).
It was further observed in \cite {\bf 10} that $B$ has real rank zero if and only if for each such $\delta$ and $x$, $vE_\delta$ can be extended to a self--adjoint unitary $u$ in $B$.
The proof of case (iii) can now proceed analogously to that of case (ii).

Thus we start with $0<\delta<1$ and $x$ in $M(A)_{\sa}$ and recursively define $w_n$ and $e_n$ in $M(I_n)_{\sa}$ and $y_n$ in $M(A)_{\sa}$ such that, with the same notations as in case (ii):
\medskip

\noindent ($12'$)\quad $w_n$ is unitary and $0\leq e_n\leq 1$,
\smallskip

\noindent ($13'$)\quad  $w_n$ extends $v_n E_\delta^{(n)}$, referring to a polar decomposition
of $\rho_n(x)$,
\smallskip

\noindent ($14'$)\quad $(\rho_{n-1} (w_n)-w_{n-1}) e_{n-1}=0$ for $n>1$,
\smallskip

\noindent ($15'$)\quad $w_n e_n=e_n w_n$
\smallskip

\noindent ($16'$)\quad $\|(\bold1-e_n) g_n\| < {1\over n}$
\smallskip

\noindent ($17'$)\quad $\rho_{n-1} (e_n) e_{n-1}=e_{n-1}$ for $n>1$
\smallskip

\noindent ($18'$)\quad $y_n\in x_{[\theta_n]}+I_n$,
\smallskip

\noindent ($19'$)\quad $\rho_n (y_n)=w_n |\rho_n (y_n)|= |\rho_n (y_n)| w_n$,
\smallskip

\noindent ($20'$)\quad $(\bold1-|y_n|) E_\delta=0$, and
\smallskip

\noindent ($21'$)\quad $(\bold1-\rho_n (|y_n|)) e_n=0$.
\medskip

To start the construction, we choose two projections $p$ and $q$ in $M(I_1)$ such that $E_{\delta\over 4}^{(1)+}\leq p\leq E_{\delta\over 8}^{(1)+}$ and $E_{\delta\over 4}^{(1)-}\leq q\leq E_{\delta\over 8}^{(1)-}$.
(The existence of $p$ and $q$ follows from the stated interpolation property and functional calculus.)
Note that $pq=0$, $\rho_1 (h_0(x_+))\in pM(I_1) p$, and $\rho_1 (h_0(x_-))\in qM(I_1) q$.
Let $(r_i)$ be an approximate identity for $p I_1 p$, $(s_j)$ an approximate identity for $q I_1 q$, $(t_k)$ an approximate identity consisting of projections for $(\bold1-p-q) I_1 (\bold1-p-q)$, and 
$$
e'_{i,j,k}=h(\rho_1 (h_0 (|x|) + (\bold1-h_0 (|x|))^{1\over 2} (r_i+s_j) (\bold1-h_0(|x|))^{1\over 2})) +t_k\,.
$$
Then $(e'_{i,j,k})$ converges strictly to $\bold1$ in $M(I_1)$.
Hence we may define $e_1=e'_{i_0,j_0,k_0}$, where $i_0,j_0,k_0$ are chosen large  enough that $\|(\bold1-e_1)g_1\|<1$.
Now let 
$$
c_1=k(\rho_1(h_0 (|x|)+(\bold1-h_0 (|x|))^{1\over 2} (r_{i_0}+s_{j_0})(\bold1-h_0 (|x|))^{1\over 2}))+t_{k_0}\,,
$$ 
$w_1=2(p+t_{k_0})-\bold1$, and $m_1=w_1c_1=c_1w_1$.
As above, $m_1=\rho_1 (x_{[\theta_1]})+z_1$ for some $z_1$ in $I_1$, and we define $y_1=x_{[\theta_1]}+z_1$.
All the conditions follow as in the proof of case (ii), or more easily.

As in the proof of case (ii), the recursive step is essentially the same as the initial step, using $y_{n-1}$ instead of $x$, $h$ instead of $h_0$, ${1\over n}$ instead of 1, and ${1\over 4}$ instead of ${\delta\over 4}$.
All the conditions are verified as in case (ii) or more easily.
And the completion of the proof after the recursive construction is also the same as in case (ii).
\medskip

(i) Let $p_n=\rho_n (p)$, $q_n=\rho_n(q)$, $B_n=p_n I_n p_n$, $C_n=q_n I_n q_n$, and $X_n=p_n I_n q_n$.
Fix strictly positive elements $g'$ of $pAp$ and $g''$ of $qAq$ and choose $g'_n$ in $B_n$ and $g''_n$ in $C_n$ as in 3.12.
Choose $u_n$ in $M(I_n)$ such that $u_n u_n^*=p_n$ and $u_n^* u_n=q_n$, and let $\wtB_n=B_n+\Bbb C p_n$, $\wtC_n=C_n+\Bbb C q_n$, and $\wtX_n=X_n+\Bbb C u_n$.
Thus $\wtB_n, \wtC_n, \wtX_n\subset M(I_n)$, and $\wtX_n$ is a $\wtB_n-\wtC_n$ Hilbert $C^*$--bimodule.

An element $u$ of $\wtX_n$ will be called unitary if $uu^*=p_n$ and $u^* u=q_n$.
Since $\wtB_n$ is unital of stable rank one, and since the map $x\mapsto xu_n^*$ is an isomorphism of $\wtX_n$ with $\wtB_n$, the result of R\o rdam \cite {\bf 39} applies:\ If $x\in\wtX_n$ with polar decomposition $x=v|x|$, and if $\delta>0$, then there is a unitary $u$ in $\wtX_n$ such that $u E_\delta=v E_\delta$.

We shall recursively construct $x_n$ in $X_n$ and $w_n$ in $\tilde X_n$ such that:
\medskip

\noindent ($12''$)\quad $w_n$ is unitary and $\|x_n\|\leq 1$,
\smallskip

\noindent ($14''$)\quad $(w_n-w_{n-1}) |x_{n-1}|=0=|x_{n-1}^*| (w_n-w_{n-1})$ for $n>1$,
\smallskip

\noindent ($15''$)\quad $x_n=w_n |x_n|$, 
\smallskip

\noindent ($16''$)\quad $\|(\bold1-x_n^* x_n) g''_n\|$, $\|g'_n (\bold1-x_n x_n^*)\| < {1\over n}$,
\smallskip

\noindent ($17''$)\quad $(\bold1-|x_n|) |x_{n-1}|=0=|x_{n-1}^*| (\bold1-|x_n^*|)$ for $n>1$, and
\smallskip

\noindent ($21''$)\quad there is $c_n$ in $C_n$ such that $0\leq c_n\leq \bold1$ and $(\bold1-c_n) |x_n|=0$.
\medskip

To start the construction, let $w_1=u_1$ and choose an approximate identity $(r_i)$ for $C_1$.
Since $(h(r_i))$ is also an approximate identity, we may define $|x_1|=h (r_{i_0})^{1\over 2}$, where $i_0$ is chosen large enough that $\|(\bold1-h (r_{i_0})) g''_1\|$, $\|g'_1 (\bold1-w_1 h(r_{i_0}) w_1^*)\| < 1$.
Then using ($15''$) as the definition of $x_1$ and taking $c_1=k(r_{i_0})$, we have all the conditions.

For the recursive step, we take $w_n$ in $\wtX_n$ to be a unitary extension of $v E_\delta$ for some $\delta<1$, referring to a polar decomposition of $w_{n-1} c_{n-1}$.
If $(r_i)$ is an approximate identity for $C_n$, then $(h(c_{n-1}+(\bold1-c_{n-1})^{1\over 2} r_i (\bold1-c_{n-1})^{1\over 2}))$ is also an approximate identity.
Thus we may define $|x_n|=(h( c_{n-1} + (\bold1-c_{n-1})^{1\over 2}r_{i_0} (\bold1-c_{n-1})^{1\over 2}))^{1\over 2}$, where $i_0$ is chosen large enough that $\|(\bold1-|x_n|^2) g''_n\|$, $\|g'_n (\bold1-w_n |x_n|^2 w_n^*\| < 1$.
The rest is similar to the above or clear.

Then we can see, as in the proof of case
(ii), that $(x_n)$ converges strictly to an element $w$ of $M(A)$ such that $ww^*=p$ and $w^* w=q$ (also $w_n\to w)$.
\hfill $\square$
\enddemo

\bigskip

\example{3.14. Remarks} (i) We don't know whether $M(I_n)$ extremally
rich for all $n$ implies $M(A)$ extremally rich.  Nevertheless 3.13 can
be used in conjunction with 4.8 below to help prove extremal richness
for some multiplier algebras, as in 4.9.
\medskip

\noindent (ii) Case (i) of 3.13 is relevant in the context of low
rank because there are known relationships, and interest in
investigating possible further relationships, between low rank and
various cancellation properies for equivalence classes of
projections.
\endexample

\vskip2truecm

\subhead{4. PULLBACKS, EXTENSIONS AND LOW RANK}\endsubhead

\bigskip

Recall that the meaning of {\it pullback diagram}, in the notation
below, is that $(\eta ,\rho )$ gives an isomorphism of $A$ with
$B\oplus_D C=\{(b,c)\in B\oplus C\mid \tau (b)=\pi (c)\}$.

\bigskip

\proclaim{4.1. Theorem} Consider a pullback diagram of 
$C^*$-algebras
$$
\CD
A @>>{\eta}> B\\
@VV{\rho}V  @VV{\tau}V\\
C @>{\pi}>> D
\endCD
$$
in which $\pi$ (hence also $\eta$) is surjective. Then:
\medskip

\itemitem {\rm (i)} $\tsr(A)\le \max(\tsr(B) ,\tsr(C) ))$.
\medskip

\itemitem {\rm (ii)} If $B$ and $C$ have real rank zero, then $A$ has real rank zero.
\medskip

\itemitem {\rm (iii)} If $B$ and $C$ (hence also $D$) are extremally rich, then:
\smallskip

\qquad $A$ is extremally rich and $\rho$ extreme-point-preserving (e.p.p.)$\Leftrightarrow \tau$ is e.p.p.
\endproclaim

\demo{Proof} By forced unitization we may assume that 
all $C^*$-algebras and all morphisms are unital. 

Note first that $\rho|\ker\eta$ is an isomorphism onto 
$\ker\pi$, so that we may put $I = \ker\eta = \ker\pi$ 
to obtain the commutative diagram of extensions 
$$
\CD
0 @>>> I @>>> A @>>{\eta}> B @>>> 0\\
@.     @|     @VV{\rho}V   @VV{\tau}V\\
0 @>>> I @>>> C @>{\pi}>> D @>>> 0
\endCD
$$
cf\. \cite{{\bf 35}, Remark 3.2}.

(i) Let $d=\max (\tsr(B) ,\tsr(C) )$.  If $\underline x$ is a tuple in
$A^d$, then we first approximate $\eta (\underline x)$ by a
unimodular tuple in $B^d$.  By the definition of quotient norm, we
approximate $\underline x$ by a tuple $\underline y$ such that
$\eta (\underline y)$ is unimodular.  Then 
$\pi (\rho (\underline y))=\tau (\eta (\underline y))$, which is
unimodular.  By \cite{{\bf 23}, Lemma 2.1} we can approximate
$\rho (\underline y)$ by a unimodular tuple $\underline z$ such that
$\pi (\underline z)=\pi (\rho (\underline y))$.  Then the pair
$(\eta (\underline y), \underline z)$ satisfies the pullback
condition and gives a unimodular approximant to $\underline x$.

(iii) If $\tau$ is e.p.p. and $u\in \Cal E(B)$, then since $C$ is
extremally rich there is $v$ in $\Cal E(C)$ such that
$\pi (v)=\tau (u)$.  Then the pair $(u,v)$ represents a lifting of
$u$ to $\Cal E(A)$.

To show that $\rho$ is e.p.p.,
consider $w=(u,v)$ in $\Cal E(A)$. Then $v$ 
is a partial isometry in $C$ such that 
$(\bold1-v^*v)C(\bold1-vv^*) \subset I$, since 
$\pi(v)=\tau(u)\in \Cal E(D)$. However, 
$(\bold1-v^*v)I (\bold1-vv^*)= 0$ because $I\subset A$ 
and $w\in\Cal E(A)$. Taken together this means that 
$v$ is in $\Cal E(C)$ as desired.

According to \cite {{\bf 6}, Theorem 6.1}, to finish the proof
that $A$ is extremally rich, we must 
check that 
$$
I + \Cal E(A) \subset (A_q^{-1})^=\,.
$$
For this, consider $w=(u,v)$ in $\Cal E(A)$ and $x$ in 
$I$. Since $C$ is extremally rich $v+x\in (C^{-1}_q)^=$, 
so $v+x$ is the limit of a sequence $(a_n)$ from $C^{-1}_q$. 
By \cite {{\bf 9}, 2.13} we may assume that $v-a_n\in I$. In the
standard decomposition $a_n = v_n e_n$, with $v_n$ in 
$\Cal E(C)$ and $e_n=(|a_n |+\bold1-v_n^*v_n)\in C^{-1}_+$ (cf\. \cite{{\bf 6}, 
Theorem 1.1}), we then have $\pi(e_n)=\bold1$. It 
follows that $w_n = (u,v_n)\in \Cal E(A)$ and $x_n =
(\bold 1,e_n) \in A^{-1}_+$; and since $w_nx_n\to(u,v+x)$ 
we have shown that $(u,v+x)\in (\QA)^=$, as desired. 

Finally, to show the reverse implication, assume 
$u\in \Cal E(B)$.  Then $u=\eta (w)$, $w\in \Cal E(A)$, since $A$ is
extremally rich.  Then by hypothesis, $\rho (w)\in \Cal E(C)$, and
hence $\tau (u)=\pi (\rho (w))\in \Cal E(D)$.

(ii) The proof is similar to, and slightly easier than, the first
part of the proof of (iii).  By \cite{{\bf 5}, Theorem 3.14} we need 
only prove that projections lift from $B$ to $A$.
\hfill{$\square$} 
\enddemo
 
\bigskip

\example{4.2. Remarks and Example} A number of papers  
contain results similar 
to ours. Thus \cite{{\bf 40}, Corollary 3.16} and
\cite{{\bf 23}, Corollary 2.7} cover Theorem 4.1 for stable rank, and 
\cite{{\bf 25}, Lemma 1.3} covers it for real rank, in the special 
case where both $\pi$ and $\tau$ (hence also $\eta$ 
and $\rho$) are surjective. Sheu's, Nistor's, and Osaka's results 
cover arbitrary values of the rank, not just low ranks.  

More recently, independently of our result (and after it was 
first obtained in 1998), Nagisa, Osaka, and Phillips \cite{{\bf 21}, Proposition 1.6}
proved 4.1(ii) for arbitrary values of the real rank.  After \cite{{\bf 21}, Corollary 1.12} they remark that their proof (which is considerably longer than ours)
also works for stable rank.  The idea, used in 4.4 below, of combining 4.1
with Busby's  analysis of extensions, is also found in Osaka's survey article
\cite{{\bf 26}, Proposition 3.4}.

The surjectivity condition cannot be entirely dropped from
Theorem 4.1.
To show that surjectivity
of $\pi$ cannot be omitted for the case of real 
rank zero take $A = C([0,1])$, embedded in the 
algebra $D$ of all bounded functions on [0,1]. Let
$C$ be the subalgebra of $D$ consisting of functions 
that are continuous on [0,1], except for possible 
jump discontinuities at points of the form
$n2^{-m}$, where $0 < n < 2^m, \,m\in \Bbb N$. Let $B$ 
be defined as $C$, except that the jump 
discontinuities are now allowed at points of the 
form $n3^{-m}$, where $ 0 < n < 3^m,\, m\in\Bbb N$. 
Realizing $B$ and $C$ as inductive limits of 
algebras of step-functions with only finitely many 
jumps, we see easily that they both have real rank 
zero. In fact, $B$ and $C$ are both isomorphic to 
$C(\Cal C)$, where $\Cal C$ denotes the Cantor set.  
Since $B\subset D$ and $C\subset D$ it is easy to 
verify that we have $B\oplus_D C 
= B\cap C = A$. But $A$ has real rank one, not zero.   

Tensoring $A,B,C$ and $D$ with $C([0,1])$, we obtain 
an example for functions on the unit square, which 
shows that surjectivity of $\pi$ can not be left
out in the stable rank one or extremal richness cases either. 
   
If $\tau$ in 4.1(iii) is not e.p.p., it can actually happen that
$\rho$ is e.p.p. and $A$ not extremally rich (cf\. 4.4) or that
$A$ is extremally rich and $\rho$ not e.p.p.
\endexample

\bigskip 

\proclaim{4.3. Corollary} Let $I$ be an ideal in 
an extremally rich $C^*-$algebra $A$ and denote by 
$\pi\colon A \to A/I$ the quotient morphism. Then for 
each extremally rich $C^*-$subalgebra $B$ which is
e.p.p. embedded in $A/I$, the $C^*-$subalgebra 
$\pi^{-1}(B)$ is extremally rich and e.p.p. embedded in $A$.
\endproclaim

\demo{Proof} We have the commutative diagram
$$
\CD
I @>>> \pi^{-1}(B) @>>{\pi}> B \\
@|     @VVV   @VVV\\
I @>>> A @>{\pi}>> A/I
\endCD
$$
which shows that $\pi^{-1}(B)$ is one of the pullbacks 
covered by Theorem 4.1.
\hfill $\square$
\enddemo

\bigskip

Every extension of $C^*-$algebras is associated with 
a Busby diagram, cf\. \cite{{\bf 12}} or \cite{{\bf 14}}, 
$$
\CD
0 @>>> I @>>> A @>>{\eta}> B @>>> 0\\
@.     @|     @VV{\rho}V   @VV{\tau}V\\
0 @>>> I @>>> M(I) @>{\pi}>> Q(I) @>>> 0
\endCD
$$
where $Q(I)=M(I)/I$ denotes the corona algebra of $I$. 
Here the right hand square is a pullback, and $A$ is 
completely determined by the Busby invariant $\tau$. 
Either $A$ is unital, which implies that also 
$\tau$ is unital; or $A$ is non-unital, 
in which case we obtain a new pullback diagram 
replacing $\eta$ and $\tau$ by the forced unitized
 morphisms $\widetilde\eta\colon \wtA\to\wtB$ and 
$\widetilde\tau\colon\wtB\to Q(I)$. Since this will 
not effect the rank of any of the algebras involved, 
we may as well assume that the extension is unital.

Applying Theorem 4.1 to the pullback diagrams described 
above, we obtain a simple but powerful tool for 
producing examples of extremally rich $C^*-$algebras.  

\bigskip 

\proclaim{4.4. Corollary} Consider an extension of 
$C^*-$algebras
$$
0 @>>> I @>>> A @>>> B @>>> 0 
$$
determined by the Busby invariant $\tau\colon B\to Q(I)$. Then
\medskip

\itemitem {\rm (a)} $\tsr(A) \le \max(\tsr(B) ,\tsr(M(I)))$.
\medskip

\itemitem {\rm (b)} If both $B$ and $M(I)$ have real rank zero, then $A$
has real rank zero.
\medskip

\itemitem {\rm (c)} If both $B$ and $M(I)$ are extremally rich, then the
following are equivalent:
\smallskip

\qquad \,\,\,{\rm (i)}\quad $\tau$ is e.p.p.
\smallskip

\qquad \,\,{\rm (ii)}\quad $\eta (\Cal E(A))=\Cal E(B)$, where $\eta\colon A \to B$ is the quotient map.
\smallskip

\qquad {\rm (iii)}\quad $A$ is extremally rich.
\endproclaim

\demo{Proof} In this situation the map $\rho$ is always e.p.p.
In fact if $u\in \Cal E(A)$, then 
$(\bold1 -u^*u)I(\bold1 -uu^*)=0$.  Since $\rho (I)$ is strictly 
dense in $M(I)$, this implies that 
$(\bold1 -\rho (u)^*\rho (u))M(I)(\bold1 -\rho (u)\rho (u)^*)=0$.
The rest follows from 4.1 and its proof.
\hfill$\square$
\enddemo

\bigskip

Sometimes it is useful to know whether
$\eta (\Cal E(A))=\Cal E(B)$ for purposes other than determining
whether $A$ is extremally rich.  For example, we may want to
know whether $\eta$ induces a surjective map on the extremal
$K-$sets, cf\. \cite {\bf 8}.  Also if $I$ is the largest ideal
of $A$ which is a dual $C^*-$algebra, then the description of
elements of $A$ with persistently closed range, cf\.
\cite{{\bf 8}, \S 7}, can be simplified if extremals lift.  Thus
it is worthwhile to state the following corollary, which follows from
the arguments above.

\bigskip

\proclaim{4.5. Corollary} If $I$ is an ideal in a unital $C^*-$algebra
$A$, $\eta\colon A \to A/I$ is the quotient map,
and if $M(I)$ is extremally rich, then $\eta (\Cal E(A))=\Cal E(A/I)$
if and only if the Busby invariant $\tau\colon A/I\to Q(I)$ is 
extreme-point-preserving.
\endproclaim

\bigskip

\example{4.6. Remarks} (i) For stable rank one or real rank zero the 
hypotheses of 4.4 can actually be weakened: If $\tsr (B)=\tsr(I)=1$ 
and the map $\iota_0\colon K_0(I)@>>>K_0(M(I))$ is injective 
then $\tsr(A)=1$. It is well-known by now that the 
vanishing of the map $\partial_1\colon K_1(B)@>>> K_0(I)$
suffices for an extension of stable rank one algebras 
to have stable rank one, and this fact follows from 
the injectivity of $\iota_0$ by an easy diagram chase. 
According to \cite{\bf 24} 
the statement about $\partial_1$ is an unpublished 
result of G. Nagy, who later published a different proof in \cite{{\bf 22}, Corollary 2}. 
Similarly, if $\RR(B)=0=\RR(I)$ and the map 
$\iota_1\colon K_1(I)@>>> K_1(M(I))$ is injective then 
$\RR(A)=0$. Here we need the injectivity of $\iota_1$ to 
show that $\partial_0 \colon K_0(B)@>>> K_1(I)$ vanishes. 
This statement is due to Zhang, cf\. \cite{{\bf 5}, 
Propositions 3.14 \& 3.15}.
\medskip

\noindent (ii) Since any morphism from an 
isometrically rich $C^*-$algebra is extreme-point-preserving, 
Corollary 4.4(c) yields extremal richness for $A$
whenever $B$ is isometrically rich, 
and this is not a severe restraint. By contrast, the 
condition that not only $I$ but also $M(I)$ 
should be extremally rich or of low rank is quite 
restrictive.
\endexample

\bigskip 

\example{4.7. Examples} To apply Corollary 4.4(c) 
we need a supply of $C^*-$algebras 
(necessarily extremally rich) whose multiplier algebras
are extremally rich. As shown in \cite{{\bf 18}, 
Corollary 3.8} this happens for every $\sigma-$unital 
purely infinite simple $C^*-$algebra. However, if $A$
is $\sigma-$unital, simple (but not elementary) and has a 
finite trace then neither $M(A)$ nor $M(A\otimes\Bbb K)$ 
are extremally rich by \cite {{\bf 18}, Theorems 3.1 \& 3.2}. 
Secondly, we immediately observe that if 
$A$ is a {\it dual} $C^*-$algebra, then $M(A)$ is extremally 
rich. Indeed, $A$ is the direct sum $\bigoplus A_i$ of 
elementary $C^*-$algebras (full matrix algebras or algebras 
of compact operators on some Hilbert space), so $M(A)$ is 
the direct product $\prod M(A_i)$, where each $M(A_i) = 
\Bbb B(\Cal H_i)$ for some Hilbert space $\Cal H_i$.  But 
if $A = c\otimes \Bbb K$, the algebra 
of norm convergent sequences of compact operators on 
$\ell^2$, then $M(A)$ is not extremally rich. This will 
be shown in Example 4.10, below.
Finally, Theorem 5.9 below provides additional examples.
\endexample

\bigskip

\proclaim{4.8. Theorem} Let $I$ be an ideal in a 
$\sigma-$unital $C^*-$algebra $A$. Then:
\medskip

\itemitem {\rm (i)} If $\RR (M(I))=0$ and $\RR (M(A/I))=0$, 
then $\RR(M(A))=0$.
\medskip

\itemitem {\rm (ii)} If $\tsr (M(I))=1$ and $\tsr (M(A/I))=1$, 
then $\tsr(M(A))=1$.
\medskip

\itemitem {\rm (iii)} If both $M(I)$ and $M(A/I)$ are 
extremally rich and $\tsr(M(A/(I+I^\perp)))=1$, then 
$M(A)$ is extremally rich and the natural morphism 
$\overline\rho\colon M(A)@>>> M(I)$ is 
extreme-point-preserving.
\endproclaim

\demo{Proof} Setting $B=A/I$ we have a commutative 
diagram in which each row contains an extension and 
$\iota$ denotes an unspecified embedding:
$$
\CD
0 @>>> I @>{\iota}>>   A   @>{\eta}>>               B @>>> 0 @.     {} \\
@. @|         @VV{\rho}V                  @VV{\tau}V  @.           @. \\
0 @>>>I @>{\iota}>>  \rho(A) @>{\pi}>>           \tau(B) @>>> 0    @.  {} \\
@. @|         @VV{\iota}V                   @VV{\iota}V    @.    @. \\
0 @>>>I @>{\iota}>>  M(I) @>{\pi}>>              Q(I) @>>> 0    @. {} \\
@. @VV{\iota}V    @AA{\iota}A               @AA{\iota}A   @.   {}      @. \\
0 @>>> M(\rho(A),I) @>{\iota}>> I(\rho(A)) @>{\pi}>> I(\tau(B)) @>{\phi}>>
M(\tau(B)) @>>> 0 \\
@. @AA{\iota}A  @AA{\overline\rho}A @.  @AA{\overline\tau}A @. \\ 
0 @>>> M(A,I) @>{\iota}>> M(A) @.
                                 \overset{\overline\eta}\to
  {\hbox to 0pt{\hss\hbox to 10.4em{\rightarrowfill}\hss}}
                                                    @. M(B)  @>>> 0
\endCD
$$
Here $M(A,I)=\{x\in M(A)\mid xA+Ax\subset I\}$ and 
$M(\rho(A),I)=\{x\in M(I)\mid x\rho(A)+\rho(A)x\subset I\}$, 
whereas $I(\rho(A))$ and $I(\tau(B))$ denote the 
idealizers of the two algebras inside $M(I)$ and $Q(I)$, 
respectively. The two quotient morphisms 
$\eta$ and $\pi$ are central to the picture and describe the 
extensions in the first and third row. The morphism 
$\rho\colon A@>>> M(I)$ is the natural map arising from the 
embedding of $I$ as an ideal in $A$; and $\ker\rho = I^\perp$. 
The Busby invariant $\tau$ is derived from $\rho$ to make the 
upper right square commutative. The morphisms $\overline\eta, 
\overline\tau$ and $\overline\rho$ are the canonical extensions 
to the multiplier algebras of $\eta, \tau$ and $\rho$, 
respectively. Since $A$, hence also $B$ and $\tau(B)$ are 
$\sigma-$unital, the overlined morphisms are surjective by 
\cite{{\bf 32}, Theorem 10}. Finally, the natural morphism $\phi$
is surjective by \cite{{\bf 14}, Corollary 3.2}. 

The kernel of the morphism $\overline\rho$ is
$$
\ker\overline\rho = \{x\in M(A)\mid xA+Ax \subset I^\perp\} = 
M(A,I^\perp)\,,
$$
which intersects $M(A,I)$ in $0$. Thus $\overline\rho$ gives an 
isomorphism of $M(A,I)$ onto the hereditary $C^*-$subalgebra 
$M(\rho(A),I)$ of $M(I)$, which is an ideal of $I(\rho(A))$. 
In the diagram we may therefore identify the two isomorphic 
ideals in $M(A)$ and $I(\rho(A))$.

We claim that the lower right rectangle is a pullback, so that 
$$
M(A)= I(\rho(A))\oplus_{M(\tau(B))}M(B)\,.
$$
For this it suffices to show that $M(\rho(A),I)= 
\ker(\phi\circ\pi)$ in $I(\rho(A))$, cf\. \cite{{\bf 35}, 
Proposition 3.1}. But this is evident, since the kernel of 
$\phi$ is the (two-sided) annihilator $\tau(B)^\perp$ in 
$I(\tau(B))$ 
\medskip

\noindent (i). Since $M(\rho(A),I)$ is hereditary
in $M(I)$, 
it follows that $\RR(M(A,I))=0$.
Given that also $\RR (M(B))=0$ we need only show that 
projections lift from $M(B)$ to $M(A)$, cf\. \cite{{\bf 5}, 3.14}.
Given a projection $p$ in $M(B)$ let $\overline h$ be a 
self-adjoint lift in $M(A)$ of the symmetry $2p-1$, and put 
$h=\overline\rho(\overline h)$. By \cite{{\bf 4}}, see also 
the explanatory version in \cite{{\bf 10}}, there is an 
interpolating projection  $q$ in $M(I)$ such that if $f_{\pm}$
are (any) two continuous functions vanishing on $[-1, \tfrac 12]$
and $[-\tfrac 12, 1]$, respectively, with $f_+(1)=1$ and
$f_-(-1)=1$ for future use, then
$$
qf_+(h)=f_+(h)\quad\text{and}\quad (\bold1-q)f_-(h)=f_-(h).
$$
For every $b$ in $B$ we have
$$
\pi(f_+(h))\tau(pb)=\tau(pb)\quad\text{and}\quad
\pi(f_-(h))\tau((\bold1-p)b)=\tau ((\bold1-p)b)\,.
$$
Consequently,
$$
\pi(q)\tau(pb)=\pi(qf_+(h))\tau(pb)=\tau(pb)\,.
$$
Similarly $\pi(\bold1-q)\tau((\bold1-p)b)=\tau((\bold1-p)b)$. Taken 
together, this means that $\pi(q)\tau (b)=\tau(pb)$.
Similarly $\tau(b)\pi(q)=\tau(bp)$. Thus $\pi(q)\in 
I(\tau(B))$ and $\phi(\pi(q))=\overline\tau(p)$. It 
follows that $\overline p=(q,p)\in M(A)$ and is a 
projection lift of $p$, as desired.
\medskip

\noindent (ii). We now know that both $M(B)$ and 
$M(A,I)$ (being isomorphic to a hereditary 
$C^*-$subalgebra of $M(I)$) have stable rank one, so 
in order to prove that $\tsr M(A)=1$ we need only 
show that unitaries lift from $M(B)$ to $M(A)$, cf\. 
\cite{{\bf 6}, 6.4}.
Given a unitary $u$ in $M(B)$ let $\overline h$ be 
a lift of $u$ to $M(A)$ and put $h=\overline\rho
(\overline h)$. Since $\tsr M(I)=1$ there is, by 
\cite{{\bf 39}, Theorem 2.2} or \cite{{\bf 33}, 
Corollary 8}, for any continuous function $f$ 
vanishing on $[0, \tfrac 12]$, with $f(1)=1$, a 
unitary $w$ in $M(I)$ such that if $h=v|h|$ is 
the polar decomposition in some $B(H)$ (cf\. the proof of 3.13) then 
$$
vf(|h|)=wf(|h|) \quad\text{and}
\quad f(|h^*|)v=f(|h^*|)w\,. 
$$

\noindent For each $b$ in $B$ we compute
$$
\pi(f(|h|))\tau(b)=\tau(f(|u|)b)=\tau (b)\,.
$$
Similarly $\tau(b)\pi(f(|h^*|))=\tau(b)$. Consequently,
$$
\pi(w)\tau(b)=\pi(wf(|h|))\tau(b)=\pi(vf(|h|))\tau(b)
=\tau (ub)\,,
$$ 
and similarly $\tau(b)\pi(w)=\tau(bu)$. Thus $\phi(\pi(w))=
\overline\tau(u)$ and 
$\overline u =(w,u)\in M(A)$ and is a unitary lift of $u$, as desired.
\medskip

\noindent (iii). Now $M(A,I)$ and $M(B)$ are 
both extremally rich, so to prove that  $M(A)$ is 
extremally rich we must show that extreme partial 
isometries in $\Cal E(M(B))$ lift in a ``good'' way,
cf\. \cite{{\bf 6}, Theorem 6.1}.  
Given $u$ in $\Cal E(M(B))$ let $\overline h$ be a lift in 
$M(A)$ and put $h=\overline\rho(\overline h)$ as in case (ii). 
Since $M(I)$ is extremally rich there is for any continuous 
function $f$ vanishing on $[0, \tfrac 12]$, with $f(1)=1$,
an extreme partial isometry  $w$ in $\Cal E(M(I))$ such 
that if $h=v|h|$ is the polar decomposition, 
then 
$$
vf(|h|)=wf(|h|)\quad\text{and}
\quad f(|h^*|)v=f(|h^*|)w\,,
$$
cf\. \cite{{\bf 6}, Theorem 2.2}.
We have assumed that $\tsr(M(\tau(B)))=1$, so $\overline\tau(u)$ 
is unitary. Consequently 
$$
\pi(f(|h|))\tau(b)=\tau(f(|u|)b)=\tau (b)
$$
for every $b$ in $B$, and similarly 
$\tau(b)\pi(f(|h^*|))=\tau(b)$. As in case (ii) we 
therefore obtain the equations 
$$
\pi(w)\tau(b)=\tau(ub)\quad\text{and}
\quad \tau(b)\pi(w)=\tau(bu)\,,
$$
so that $\pi(w)\in I(\tau(B))$ with $\phi(\pi(w))=
\overline\tau(u)$. Thus $\overline u = (w,u)\in M(A)$ 
and is an extremal lift of $u$ as we wanted.

Let $p_{\pm}$ denote the defect projections of 
$\overline u$. The extra technical condition needed 
in \cite{{\bf 6}, Theorem 6.1} is that the two 
bimodules $p_{\pm}M(A)r$ are extremally rich for any 
defect projection $r$ arising from an element 
$s$ in $\Cal E(M(A,I)\,\widetilde{}\,)$.
Since $\overline\tau(\overline\eta(\overline u))
=\overline\tau(u)$ is unitary it follows that 
$\pi(\overline\rho(p_{\pm}))\in \tau(B)^\perp$, 
i\.e\. $\overline\rho(p_{\pm})\in M(\rho(A),I)$.
Since moreover $r\in M(A,I)$ also $\overline\rho(r)\in 
M(\rho(A),I)$. But this is a hereditary 
$C^*-$subalgebra of $M(I)$, so 
$$
\overline\rho
(p_\pm M(A)r)=\overline\rho(p_\pm)M(\rho(A),I)\overline\rho(r)
=\overline\rho(p_\pm)M(I)\overline\rho(r)\,.
$$
Note now that by construction $\overline\rho(\overline u)
=w \in \Cal E(M(I))$, so $\overline\rho(p_\pm)$ are both 
extreme defect projections of $M(I)$. Moreover, since 
$\overline\rho(M(A,I))$ is hereditary in $M(I)$, $\overline\rho(s)$
is in $\Cal E(M(I))$ and $\overline\rho(r)$ is also an
extreme defect projection of $M(I)$.
It now follows from \cite{{\bf 6}, Proposition 4.4}
that $\overline\rho(p_\pm)M(I)\overline\rho(r)$ is 
extremally rich. 

Finally, the fact that $\overline\rho$ is e.p.p. follows by
the same argument as in 4.4.
\hfill$\square$
\enddemo

\bigskip

\proclaim{4.9. Corollary} Let $\{I_\alpha \mid 0\le \alpha \le \beta \}$
be a composition series for a separable $C^*-$algebra $A$.  Then
\medskip

\itemitem {\rm (i)} If $\RR(M(I_{\alpha +1}/I_\alpha))=0$ for all 
$\alpha <\beta$, then $\RR(M(A))=0$.
\medskip

\itemitem {\rm (ii)} If $\tsr(M(I_{\alpha +1}/I_\alpha))=1$ for all
$\alpha <\beta$, then $\tsr(M(A))=1$.
\medskip

\itemitem {\rm (iii)} If $M(I_1)$ is extremally rich and $\tsr(M(I_{\alpha +1}/I_\alpha))=1$ for $1\le \alpha <\beta$, then $M(A)$ is extremally rich.
\endproclaim

\demo{Proof} If we assume, as we may, that $(I_\alpha)$ is strictly 
increasing, then $\beta$ is countable, since $\hA$ is second
countable.  Then the result follows by a routine transfinite 
induction from Theorems 3.13 and 4.8.
\hfill $\square$
\enddemo

\bigskip

\example{4.10. Examples} In \cite{{\bf 8}, Example 7.9} 
we explored a non-extremally rich $C^*-$algebra $A$ that illustrates a 
number of points in the present paper.  (There is a typographical
error in the second-to-last paragraph of \cite{{\bf 8}, Example 7.9}: 
Please change (vii) to (viii).)  Let $I=
\bigoplus_{n=1}^\infty \Bbb M_n$ and let $B$ be the 
$C^*-$subalgebra of $M(I)\,(=\prod_{n=1}^\infty \Bbb M_n)$ 
consisting of norm convergent 
sequences (relative to the standard embedding of 
$\Bbb M_n$ in $\Bbb B(\ell^2))$. Then $M(B)$ 
consists of the strong*  convergent sequences. If 
$s_n$ denotes the forward truncated shift in $\Bbb M_n$ and
$s$ the forward shift on $\ell^2$, we put 
$\widetilde s = (s_n)$ in $M(B)$ and now $A=
C^* (B\,, \widetilde s\,, \bold 1)$ in $M(B)$. The quotient $A/I$ 
is isomorphic to $\Cal T_e$, the extended Toeplitz
algebra generated by the element  $s\oplus s^*$ in
$\Bbb B(\ell^2\oplus\ell^2)$. The 
primitive ideal space of $A$ (equal to the spectrum $\widehat A$) 
is the disjoint union 
$$
\{\pi_1, \pi_2, \dots\}\cup\{\sigma_+\}
\cup\{\sigma_-\}\cup\Bbb T\,,
$$
where each $\pi_n$ corresponds to an $n-$dimensional 
representation, $\sigma_\pm$ correspond to  the
infinite dimensional representations of $\Cal T_e$, and 
the circle $\Bbb T$ (with the usual topology) 
consists of one-dimensional representations. The set 
$\Cal F=\{\sigma_+, \sigma_-\}\cup\Bbb T$ is the hull of $I$, and
$(\pi_n)$ converges to all points in $\Cal F$ simultaneously. 

Each set $\{\pi_n\}$ is closed in $\hA$ and 
corresponds to a stable rank one quotient of
$A$. Thus the hypothesis in (4) of 2.13(ii) that $A/I_1$ be
isometrically rich cannot be weakened to extremal richness.
Also, each of the sets $\{\sigma_+\}\cup\Bbb T$ and
$\{\sigma_-\}\cup\Bbb T$ is closed in $\hA$ and
corresponds to an isometrically rich quotient of
$A$. Thus the disjointness condition in (5) of 2.13(ii) cannot 
be omitted.  Since $\widehat A$ is an 
almost Hausdorff space, the Hausdorff demand in 
part (iv) of Theorem 2.11 cannot be replaced by 
almost Hausdorff.

Since $A/B$ is isometrically rich and $A$ is not
extremally rich, it follows from 4.4 that $M(B)$ is not extremally
rich.  This shows that the non-existent parts (iii$'$) and (iv$'$) of 2.12
are false if we consider 
$B^\vee =(\bigcup_n \{ \pi_n \} )\bigcup \{\sigma_+\}$.  Note that $B^\vee$ is Hausdorff and all but one of the quotients of $B$ are
unital and of stable rank one.  This also shows that
the hypothesis that
$\tsr(M(A/(I+I^\perp)))=1$ cannot be dropped from
Theorem 4.8, since $\tsr(M(I))=1$ and $M(B/I)$ is isometrically
rich, thus suggesting some sharpness in the hypotheses of 4.8.
Finally, since $B$ is a corner of $c\otimes\Bbb K$, 
the assertion in Example 4.7 that 
$M(c\otimes\Bbb K)$ is not extremally rich has 
been verified.
\endexample

\vskip2truecm

\subhead{5.  APPLICATIONS, REMARKS, AND QUESTIONS}\endsubhead

\bigskip

The main purpose of this section is to provide some sample applications of the basic results of the paper by determining when $CCR$ algebras or
their multiplier algebras have low rank.
In deciding how much material to present and how to present it, we have tried to walk a fine line.
On the one hand we don't want to obscure the main purpose with too many
technical proofs, and on the other hand we don't
want to complicate the statements of the results with unnecessary technical hypotheses.
But we begin with two light contributions concerning type $I$ algebras.
These together with 2.9 constitute our best effort to give a somewhat general characterization of low topological dimension.
(The ``lightness'' of 5.1 lies in the fact that the most important parts were already known.)

We will use the local definition of $AF$--algebra in order to cover non--separable algebras.
Thus $A$ is an {\it $AF$--algebra} if for every finite subset $\Cal F$ of $A$ and every $\epsilon>0$, there is a finite dimensional
$C^*$--subalgebra $B$ such that dist$(a,B)<\epsilon$ for each $a$ in $\Cal F$.
In 5.1 below the equivalence of (i) and (iv) in the separable case was proved by Bratteli and Elliott \cite {\bf 2}, the fact that (i) implies (v)
is a special case of results of Lin \cite {\bf 19, 20}, and Pasnicu \cite{{\bf 28},
Remark 2.12} shows that (i) is equivalent to the ideal property in the separable case.

\bigskip

\proclaim{5.1. Proposition}Let $A$ be a type I $C^*$--algebra.
Then {\rm (v)} $\Rightarrow$ {\rm (i)} $\Leftrightarrow$ {\rm (ii)} $\Leftrightarrow$ {\rm (iii)} $\Leftrightarrow$ {\rm (iv)}.
If also $A$ is $\sigma$--unital, then all the conditions are equivalent.
\medskip

\itemitem{\rm (i)}$A$ is an $AF$--algebra.
\smallskip

\itemitem{\rm (ii)}$A$ has real rank zero.
\smallskip

\itemitem{\rm (iii)}$A$ has generalized real rank zero.
\smallskip

\itemitem{\rm (iv)}$\topdim(A)=0$.
\smallskip

\itemitem{\rm (v)}$M(A)$ has real rank zero.
\endproclaim

\demo{Proof}The implications (i) $\Rightarrow$ (ii) $\Rightarrow$ (iii) and (v) $\Rightarrow$ (ii) are obvious, and (iii) $\Rightarrow$ (iv) follows from 2.9.

(iv) $\Rightarrow$ (i):\ $A$ has a composition series $\{I_\alpha\mid 0\leq\alpha\leq\beta\}$ such that for each $\alpha<\beta$ $I_{\alpha+1}/I_\alpha$
is Rieffel--Morita equivalent to a commutative algebra $C_0(\Cal X_\alpha)$.
Each $\Cal X_\alpha$ is totally disconnected, and hence $C_0(\Cal X_\alpha)$ is $AF$.
Therefore $I_{\alpha+1}/I_\alpha$ is $AF$.
Since direct limits of $AF-$algebras are $AF$, and since an extension of one $AF-$algebra by another is $AF$ (a fact which was not known when \cite {\bf 2} was
written), a routine transfinite induction shows that $A$ is $AF$.

Finally, if $A$ is $\sigma$--unital, the implication (i) $\Rightarrow$ (v) follows from \cite {{\bf 19}, Corollary 3.7}.
\hfill{$\square$}
\enddemo

\bigskip

\proclaim{5.2. Proposition}If $A$ is a type I $C^*$--algebra, then the following conditions are equivalent:
\medskip

\itemitem{\rm (i)}$A$ has generalized stable rank one.
\smallskip

\itemitem{\rm (ii)}$A$ has a composition series $\{I_\alpha\mid 0\leq\alpha\leq\beta\}$ such that $I_{\alpha+1}/I_\alpha$ is extremally rich for each
$\alpha<\beta$.
\smallskip

\itemitem{\rm (iii)}$\topdim(A)\leq 1$.
\endproclaim

\demo{Proof}The implication (i) $\Rightarrow$ (ii) is obvious, and for (ii) $\Rightarrow$ (iii) we use a composition series such that $I_{\alpha+1}/I_\alpha$
is both extremally rich and Rieffel--Morita equivalent to a commutative $C^*$--algebra $C_0(\Cal X_\alpha)$, cf\. 2.1(vi).
Then $I_{\alpha+1}/I_\alpha$ extremally rich implies $C_0(\Cal X_\alpha)$
extremally rich which implies $\tsr(C_0(\Cal X_\alpha))=1$ which implies
$\topdim(C_0(\Cal X_\alpha))\leq 1$ which implies $ \topdim (I_{\alpha+1}/I_\alpha)\leq 1$.
Then by 2.5 $\topdim (A)\leq 1$.
The proof that (iii) implies (i) is similar.
Now $\topdim (I_{\alpha+1}/I_\alpha)\leq 1$ implies $\topdim (C_0(\Cal X_\alpha))\leq 1$ which implies $\tsr(C_0(\Cal X_\alpha))=1$ which implies
$\tsr(I_{\alpha+1}/I_\alpha)=1$.
\hfill $\square$
\enddemo

If $A$ is type $I$ and $\topdim(A)=0$, then $A$ has stable rank one since it is $AF$.
But in the one--dimensional case there are numerous examples where $A$ is extremally rich but not of stable rank one and numerous examples where $A$ is not
even extremally rich.
But, as we proceed to show, such examples cannot be $CCR$ algebras, despite the fact that the (non-canonical) composition series for a $CCR$ algebra with each
$I_{\alpha+1}/I_\alpha$ of continuous trace can be very complicated.

\bigskip

\proclaim{5.3. Lemma}Let $A$ be an $n$--homogeneous $C^*$--algebra such that $\topdim(A)\leq 1$.
\medskip

\itemitem{\rm (i)}If $A$ is $\sigma$--unital, then $\tsr(M(A))=1$.
\medskip
\itemitem{\rm (ii)}In any case, $\tsr(A)=1$.
\endproclaim

\demo{Proof}(i) By the structure theory for $n$--homogeneous algebras and the fact that $\hA$ is $\sigma$--compact, we can find a sequence $(I_m)$ of ideals
such that $\bigcup \hull (I_m)=\hA$ and each $A/I_m$ is isomorphic to $C(\Cal X_m)\otimes \Bbb M_n$ where $\Cal X_m$ is compact.
Since $\dim (\Cal X_m)\leq 1$, $A/I_m$ is unital and of stable rank one.
The conclusion follows from 2.12(ii$'$).

(ii) If $A$ is $\sigma$--unital (ii) follows from (i).
But $A$ is the direct limit of an upward directed family of $\sigma$--unital ideals.
\hfill{$\square$}
\enddemo

\bigskip

\proclaim{5.4. Lemma}Assume each irreducible representation of $A$ has dimension at most $n$ and $\topdim(A)\leq 1$.
\medskip

\itemitem{\rm (i)}Then $\tsr(A)=1$.
\medskip

\itemitem{\rm (ii)}If also $A$ is separable, then $\tsr (M(A))=1$.
\endproclaim

\demo{Proof}(i) We use induction on $n$.
The case $n=1$ is known.
For $n>1$ $A$ has an ideal $I$ which is $n$--homogeneous such that all irreducibles of $A/I$ have dimension at most $n-1$.
Thus $\tsr(I)=\tsr (A/I)=1$ by 5.3(ii) and induction.
To complete the proof we show that $\iota_0\colon K_0 (I)\rightarrow K_0(A)$ is injective, cf\. 4.6.
Write $I=(\bigcup I_\alpha)^=$ for an upward directed family of $\sigma$--unital ideals.
Consider
$$
K_0(I_\alpha) \overset (\iota_\alpha)_*\to\rightarrow K_0(I)\overset \iota_0\to\rightarrow K_0(A)\overset (\rho_\alpha)_*\to\rightarrow K_0(M(I_\alpha)),
$$
where $\rho_\alpha \colon A\to M(I_\alpha)$ is the natural map used above, cf\. 4.4.
The composite map is injective by 5.3(i).
Since every element of $K_0(I)$ is in the image of $(\iota_\alpha)_*$ for $\alpha$ sufficiently large, this implies $\iota_0$ is injective.

(ii) We use the same induction and the same $I$, but here
we use 5.3(i) and 4.8.
\hfill{$\square$}
\enddemo

\bigskip

\proclaim{5.5. Lemma}Let $A$ be a $C^*$--algebra, all of whose irreducible representations are finite dimensional, such that $\topdim(A)\leq 1$.
\medskip

\itemitem{\rm (i)}Then $\tsr(A)=1$.
\medskip

\itemitem{\rm (ii)}If also $A$ is separable, then $\tsr(M(A))=1$.
\endproclaim

\demo{Proof}There is a sequence $(I_n)$ of ideals such that all irreducibles of $A/I_n$ have dimension at most $n$ and $\bigcup \hull (I_n)=\hA$.
Thus (i) follows from 2.12(ii) and 5.4(i), and (ii) follows from 2.12(ii$'$) and 5.4(ii).
\hfill{$\square$}
\enddemo

\bigskip

\proclaim{5.6. Theorem}If $A$ is a $CCR$ $C^*$--algebra, then the following conditions are equivalent:
\medskip

\itemitem{\rm (i)}$A$ has stable rank one.
\smallskip

\itemitem{\rm (ii)}$A$ is extremally rich.
\smallskip

\itemitem{\rm (iii)}$\topdim(A)\leq 1$.
\endproclaim

\demo{Proof}We need only prove that (iii) implies (i).
We claim that $A=(\bigcup B_\alpha)^=$, for an upward directed family of hereditary $C^*$--subalgebras, such that each $B_\alpha$ has only finite dimensional
irreducibles.
One way to see this is to use the theory of the minimal dense ideal, $K(A)$, found in \cite {{\bf 31}, \S5.6}.
By \cite {{\bf 31}, 5.6.2} the hereditary $C^*$--subalgebra generated by any finite subset of $K(A)$ is contained in $K(A)$.
And for each $a$ in $K(A)$ and each irreducible $\pi$, $\pi(a)$ has finite rank.
Since $B_\alpha^\vee$ is an open subset of $\hA$, $\topdim(B_\alpha)\leq 1$.
Then the conclusion follows from 5.5(i) and the preservation of low rank by direct limits.
\hfill{$\square$}
\enddemo

\proclaim{5.7. Theorem}If $A$ is a $\sigma$--unital $CCR$ $C^*$--algebra, then $M(A)$ has stable rank one if and only if $\topdim(A)\leq 1$ and all irreducible
representations of $A$ are finite dimensional.
\endproclaim

\demo{Proof}If $A$ has an infinite dimensional irreducible, then $A$ has a quotient algebra isomorphic to $\Bbb K$.
Therefore $M(A)$ has a quotient isomorphic to $B(\Cal H)$.
But $\tsr(B(\Cal H))=\infty$.

For the other direction we write $A=(\bigcup A_\alpha)^=$ for a suitable upward directed family of separable $C^*$--subalgebras.
Since $A$ is $\sigma$--unital, we can choose the $A_\alpha$'s so that $M(A)=\bigcup M(A_\alpha)$.
(Note that $M(A_\alpha)\subset M(A)$ if $A_\alpha$ contains an approximate identity of $A$.)
Since $\tsr(A)=1$ (by 5.5(i) or 5.6), we can also arrange that $\tsr(A_\alpha)=1$ and hence $\topdim(A_\alpha)\leq 1$.
Then 5.5(ii) applies to each $A_\alpha$, and the conclusion follows.
\hfill{$\square$}
\enddemo

\bigskip

We are not applying the inverse limit theory here, but we could have used 4.9 (which is based on 3.13) instead of 2.12 in the proof of 5.5(ii), and we
could have used an easier inverse limit argument similar to the construction in 3.4 instead of 2.12 in the proof of 5.3(i).
A negative answer to the first question below would allow the possibility of more crucial applications of Theorem 3.13.

\example{5.8. Questions}
(i) If $A$ is a $\sigma$--unital (or separable) $C^*$--algebra (or type I $C^*$--algebra) such that $M(A)$ has stable rank one, does there necessarily exist a
sequence $(I_n)$ of ideals such that $\bigcup \hull (I_n)=\hA$ and each $A/I_n$ is unital?

\noindent (ii) Does there exist a separable $C^*$--algebra $A$ such that $\tsr(A)=1$ and $1<\tsr (M(A))<\infty$?

\noindent {\bf NOTE}:\ These questions are not conjectures, and nothing in this paper should be construed as a conjecture.
\endexample

\bigskip

\proclaim{5.9. Theorem}Let $I$ be an ideal of a $\sigma$--unital $C^*$--algebra $A$ such that $I$ is a dual $C^*$--algebra, $A/I$ has only finite dimensional
irreducible representations, and $\topdim (A/I)\leq 1$.
Then $M(A)$ is extremally rich.
\endproclaim

\demo{Proof}By 5.7 $M(A/I)$ has stable rank one.
Hence the result follows from 4.8(iii).
\hfill{$\square$}
\enddemo

\bigskip

\example{5.10. Remark}It can be shown that if $M(A)$ is extremally rich for a $\sigma$--unital $CCR$ algebra $A$, then $A$ satisfies the hypotheses of 5.9.
In other words:
\medskip

\noindent If $A$ is a $\sigma$--unital $CCR$ $C^*$--algebra, then $M(A)$ is extremally rich if and only if $\topdim(A)\leq 1$ and each infinite dimensional irreducible
representation of $A$ gives an isolated point of $\hA$.
\medskip

\noindent We are omitting the proof of the converse but point out that the algebra called $B$ in Example 4.10 is an instructive example.
Note that not every $C^*$--algebra satisfying the hypotheses of 5.9 is $CCR$.
\endexample

\bigskip

\proclaim{5.11. Corollary}Let $I$ be an ideal of a $C^*$--algebra $A$ such that $I$ has only finite dimensional irreducible representations  and $\topdim(I)\leq 1$.  Then $\tsr(A)=\tsr (A/I)$.
\endproclaim

\demo{Proof}If $I$ is $\sigma$--unital, this follows directly from 5.7 and 4.4.
The general case follows from this via standard techniques for reducing to the separable case.
In representing $A$ as the direct limit of separable algebras $A_\alpha$, we arrange directly that $\tsr(A_\alpha\cap I)=1$ rather than trying to arrange
directly that $\topdim (A_\alpha\cap I)\leq 1$.
(We don't know whether it is possible to control topological dimension in general constructions of this sort.)
\enddemo

\example{5.12. Remark}As mentioned above in 4.6, one can use $K-$theory
instead of multiplier algebras to study the stable rank of an extension
of one stable rank one $C^*-$algebra by another; and in this way one can
prove the stable rank one case of 5.11 without using multiplier algebras.
  The following result,
whose proof is omitted, can be used:
\medskip

\noindent Let $A$ be a $CCR$ algebra such that $\topdim (A)\leq 1$.
If $\alpha\in K_0(A)$ and $\pi_*(\alpha)=0$ in $K_0(\pi (A))$ for every irreducible representation $\pi$, then $\alpha=0$.
\endexample

\bigskip

\proclaim{5.13. Theorem}Let $I$ be an ideal of a $C^*$--algebra $A$ such that $I$ has only finite dimensional irreducible representations, $\topdim I\leq 1$,
and $A/I$ is extremally rich.
Then {\rm (i)}, {\rm (ii)}, and {\rm (iii)} below are equivalent and {\rm(v)} implies {\rm (i)}, {\rm (ii)}, and {\rm (iii)}.
If also $I$ is $\sigma$--unital then {\rm (iv)} is equivalent to {\rm (i)}, {\rm (ii)}, and {\rm (iii)}.
\medskip

\itemitem{\rm (i)}$A$ is extremally rich.
\smallskip

\itemitem{\rm (ii)}$\pi(\Cal E(\wtA))=\Cal E((A/I)^\sim)$, where $\pi$ is the quotient map.
\smallskip

\itemitem{\rm (iii)}$\tsr(A/(I+I^\perp))=1$.
\smallskip

\itemitem{\rm (iv)}The Busby invariant $\tau\colon A/I\rightarrow M(I)/I$ is extreme-point-preserving.
\smallskip

\itemitem{\rm (v)}$A/I$ is isometrically rich.
\endproclaim

\demo{Proof}We may assume $A$ unital.

It is obvious that (i) implies (ii).
Assume (ii) and let $u$ be in $\Cal E(A/(I+I^\perp))$.
Lift $u$ first to $\Cal E(A/I)$ and then to $v$ in $\Cal E(A)$.
Let $\rho$ be a representation of $A$ which is faithful, non-degenerate,
and atomic on $I$.
Hence $\ker(\rho)=I^\perp$.
Clearly $\rho(v)$ is unitary.
Therefore $v+I^\perp$ is unitary in $A/I^\perp$ and $u$ is unitary.
It now follows from \cite {{\bf 9}, 2.7} that $\tsr (A/(I+I^\perp))=1$.
If (iii) is true, then $\tsr(A/I^\perp)=1$ by 5.11.
Since $A$ is a pullback of $A/I$ and $A/I^\perp$, with both maps surjective, 4.1 implies that $A$ is extremally rich.
Thus (i), (ii), and (iii) are equivalent.

Now if $I$ is $\sigma$--unital, the facts that (i), (ii), and (iv) are equivalent and (v) implies (i) follow from 5.7, 4.4(c), and 4.6(ii).
The fact that (v) still implies (i), (ii), and (iii) in general follows via standard techniques for reducing to the separable case.
\hfill{$\square$}
\enddemo

\bigskip

\example{5.14. Remarks and Question}
(i) One should not view 5.13 as the prototype for results stating that certain $C^*$--algebras $I$ are universally good ideals from the point of view of
extremal richness of extensions.
The reason is that the equivalence of (ii) and (iii) is an artifact of the particular class of ideals being considered.
Instead, 4.4 and 4.6 suggest two reasonable prototypes:
$$
\align
& \text{(v)}\Rightarrow \text{(i)} \Leftrightarrow \text{(ii)}\,,\text{ or } \\
& \text{(v)}\Rightarrow \text{(i)} \Leftrightarrow \text{(ii)} \Leftrightarrow \text{(iv)}\,.
\endalign
$$
One may want technical hypotheses, such as the $\sigma-$unitality used
above, in order to prove the stronger version.
It can be shown that arbitrary purely infinite simple $C^*$--algebras satisfy 
the weaker version, cf\. 4.7.
\medskip

\noindent (ii) The fact that (v) implies (iii) amounts to the statement that $M(I)/I$ contains no proper isometries.
Our proof of this uses $\topdim (I)\leq 1$.
Could this hypothesis be omitted?  In other words:
\medskip

\noindent If $I$ is a $C^*-$algebra with only finite dimensional
irreducible representations, is $M(I)/I$ necessarily stably finite?
\endexample

\example{5.15. Final Remark}Theorems on higher real rank and stable rank
of $CCR$ algebras can also be proved using our main results.  
This will be done in a future paper.  In particular,
a generalization, which doesn't involve topological dimension, of Corollary 5.11 to arbitrary ranks will be given.
\endexample

\bigskip

{\it This article, together with} \cite {{\bf 6, 7, 8, 9, 11}},
{\it  is part of the documentation of a project begun in 1992 during
a visit of
the first author at the University of Copenhagen. The main results were
obtained in 1992, but revised in 1996, 1998 and 2000, and
a few additional improvements were obtained in 2006. The authors
gratefully acknowledge the support of the Danish Research Council (SNF).}

\vskip2truecm


\bigskip

\noindent Lawrence G. Brown           

\noindent Department of Mathematics  

\noindent Purdue University         

\noindent West Lafayette           

\noindent Indiana 47906, USA      

\noindent lgb \@ math.purdue.edu

\enddocument